\documentclass[preprint,10pt,sort&compress]{elsarticle}

\usepackage{subcaption}
\usepackage{amssymb}
\usepackage{amsmath}
\usepackage{graphicx}
 \usepackage{amsthm}
 \usepackage{bm}
 \usepackage{color}

\newtheorem{theorem}{\hspace{1em}~Theorem}[section]

\newtheorem{lemma}[theorem]{\hspace{1em}~Lemma}
\newtheorem{remark}{\hspace{1em}~Remark}[section]

\journal{*}

\begin{document}

\begin{frontmatter}



\title{Stabilization-free virtual element methods based on finite element interpolation}

\author[addr1]{Jikun Zhao}
\ead{jkzhao@zzu.edu.cn}
\author[addr1]{Wenhao Zhu}
\ead{zhuwenhao@gs.zzu.edu.cn}
\author[addr2]{Bei Zhang}
\ead{beizhang@haut.edu.cn}

\author[addr3,addr4]{Shipeng Mao\corref{cor1}}
\ead{maosp@lsec.cc.ac.cn}
\cortext[cor1]{Corresponding author.}

\address[addr1]{School of Mathematics and Statistics, Zhengzhou University, Zhengzhou 450001, PR China}
\address[addr2]{School of Mathematics and Statistics, Henan University of Technology, Zhengzhou 450001, PR China}
\address[addr3]{State Key Laboratory of Mathematical Sciences, Academy of Mathematics and Systems
Science, Chinese Academy of Sciences, Beijing 100190, PR China}
\address[addr4]{School of Mathematical Sciences, University of Chinese Academy of Sciences, Beijing
100049, PR China}

\begin{abstract}
 In this paper, we introduce a new framework for designing stabilization-free virtual element methods (VEMs) based on a finite element interpolation-based strategy, where we can simultaneously eliminate the stabilization terms in the discretizations of diffusion and reaction terms. The core idea is to construct a computable, polynomial-preserving, and norm-equivalent interpolation operator from the virtual element space to a finite element space on the simplicial submesh of each polytope. Leveraging the properties of this operator, we design two types of stabilization-free schemes. The first scheme requires the interpolation to preserve the polynomial consistency related to the bilinear forms, thereby maintaining both consistency and stability as in the standard VEM. The second scheme relaxes this consistency requirement. While it may not satisfy the standard polynomial consistency, the second scheme retains optimal convergence with simpler construction and fewer degrees of freedom. We construct such interpolation operators for both conforming and nonconforming virtual elements in two and three dimensions. These operators are then employed to realize stabilization-free schemes for conforming and nonconforming VEMs. Finally, numerical experiments confirm that the proposed stabilization-free VEMs achieve optimal convergence rates and exhibit significantly lower sensitivity to strong anisotropy compared to standard VEM. 
\end{abstract}


\end{frontmatter}

% \linenumbers

%\section{}
%\label{}


\section{Introduction}

\subsection{Background}

As a numerical method for solving partial differential equations on polygonal or polyhedral meshes, the virtual element method (VEM) has undergone significant development since its introduction in \cite{BeiraodaVeiga2013}.  For various conforming and nonconforming VEM formulations, we refer to \cite{Antonietti2023,Antonietti2024,Antonietti2022,Antonietti2016a,AyusodeDios2016,BeiraodaVeiga2016a,BeiraodaVeiga2016b,BeiraodaVeiga2020,BeiraodaVeiga2017a,BeiraodaVeiga2017,BeiraodaVeiga2019,Brezzi2013,Carstensen2023,Chen2018a,Chen2019,Chen2022a,Mascotto2018a,Wei2021,Zhao2016,Zhao2018a,Zhao2017,Zhao2020,Zhao2021,Zhao2023} and the references therein. One of its key advantages lies in the flexibility of space construction and mesh generation. Unlike the finite element method (FEM), the shape function space in VEM includes not only polynomials but also certain non-polynomial functions that cannot be directly computed. To formulate the discrete VEM scheme, computable projectors are therefore employed. To ensure coercivity,  additional stabilization terms with  suitable scaling factor must be incorporated into the discrete bilinear form, thereby satisfying the so-called stability condition.  As noted in \cite{BeiraodaVeiga2013},  the choice of stabilization terms generally depends on the problem and the degrees of freedom. Among them, the most common one is the dofi-dofi type, but needs to be carefully tuned/change in presence of ``awkwardly-shaped polygons".
 As also mentioned  in \cite{Mascotto2023}, the properties of the stabilization terms were introduced and motivated heuristically in the early stages, and the theoretical aspects were systematically investigated in recent years, see \cite{BeiraodaVeiga2017a,Brenner2017,Chen2018,Brenner2018}. Several types of stabilization terms were summarized in \cite{Mascotto2023}.  Nonetheless, an inappropriate choice can significantly degrade the performance of the VEM and may even compromise its accuracy, as detailed in
  \cite{Dassi2018,Mascotto2023}.  Consequently, there has been a growing interest in entirely eliminating stabilization terms, with several approaches proposed in recent works such as 
  \cite{Berrone2025,Chen2023a,Meng2022,Chen2023,Xu2023,Xu2024,Berrone2023,Borio2024,Berrone2024,Chen2024,Berrone2025a}.

 To eliminate the stabilization terms,  a key technique is to introduce  a computable projection $Q^K$ of the gradient of virtual element (VE) functions onto a  larger space (consisting of vector-valued polynomials or piecewise polynomials) on each element $K$, such that for each VE function $v_h$, we have the following norm equivalence 
\begin{equation}\label{eq:GradProjNormEqui}
C_1\|Q^K \nabla v_h\|_K\leq\|\nabla v_h\|_K\leq C_2\|Q^K\nabla v_h\|_K,
\end{equation}
where $C_1,C_2$ are generic constants independent of the element size.

The pioneering  work in this direction appeared in  \cite{Berrone2025} by Berrone et al, which introduced an $L^2$ 
  projection of the gradient onto a higher-order polynomial space on each element for the lowest-order conforming VEM of the Poisson equation in two dimensions (2D).  The degree of the projection space depends on the number of edges and the geometry of the polygon. Building on the same technique,  Berrone et al \cite{Berrone2025a} later analyzed a higher-order stabilization-free VE scheme for general second-order elliptic problems in 2D. The approach developed in \cite{Berrone2025,Berrone2025a} has been applied to construct stabilization-free VEMs for various problems involving gradient operators, such as Laplacian eigenvalue problems  \cite{Meng2022}, convection-diffusion eigenvalue problem \cite{Marcon2025} and plane elasticity \cite{Chen2023a}. It should be noted, however, that the theoretical analysis of this approach is difficult to extend to 3D VEMs or to nonconforming VEMs. Although Xu and Wriggers \cite{Xu2024}  extended the approach to the lowest-order stabilization-free VEM for linear elasticity in 3D, therein the approach is only studied numerically.  We also mention that  in \cite{Berrone2023} an improvement of this approach is proposed for the lowest-order VEM of a general second-order elliptical problem in 2D, where the projection space was replaced by the gradient space of harmonic polynomials to reduce quadrature-related computational cost. Based on this improvement, a lowest-order stabilization-free mixed VEM was presented for the 2D Poisson equation in \cite{Borio2024}, though its stability and convergence were analyzed only on quadrilateral meshes. 
  Another variant of this approach was developed in \cite{Berrone2024}, which introduced higher-order projections onto divergence-free polynomial vector-valued spaces.  There, the lowest-order stabilization-free VEM for 2D Poisson equation was analyzed on some particular classes of polygons, while  the coercivity of bilinear form in this scheme was proved only on general quadrilateral elements. 
\textit{In summary, the approach originating from  \cite{Berrone2025} and its variants face significant challenges in establishing well-posedness and convergence for resulting stabilization-free VEMs in 3D. 
Moreover, this approach seems to be not applicable to developing the nonconforming VEMs without stabilization.}

Chen, Huang, and Wei further advanced this direction in \cite{Chen2024} by proposing a unified framework for constructing stabilization-free conforming and nonconforming VEMs of arbitrary order for the Poisson equation in arbitrary dimension.
 Their key innovation lies in replacing the gradient projection space with a special H(\text{div})-conforming macro finite element space on a simplicial subdivision of each polytope, for which the norm equivalence  \eqref{eq:GradProjNormEqui}  is rigorously proved. 
 \textit{It should be noted, however, that the 3D conforming VEM presented in \cite{Chen2024}, requires all mesh faces to be triangles-a geometric constraint that ensures both the computability of VE functions on faces and the validity of the norm equivalence for the gradient projection $Q^K$.}

%Besides, Chen, Huang and Wei developed the stabilization-free conforming and nonconforming VEMs of arbitrary degree for Poisson equation in arbitrary dimension in a unified way in \cite{Chen2024}, where the gradient projection space was replaced by a special $H$(div)-conforming macro FE space on a simplicial partition of each polytope in mesh and the norm equivalence \eqref{eq:GradProjNormEqui} was rigorously proved. \textbf{\textit{However,  for the stabilization-free conforming VEM in 3D presented in \cite{Chen2024}, in order to obtain the computability of VE functions on faces and the norm equivalence of gradient projection $Q^K$, the faces in meshes must be triangles.}}

%\textbf{\textit{In addition to the aforementioned issues, the introduction of gradient projection $Q^K$ with norm equivalence \eqref{eq:GradProjNormEqui} can only be utilized to eliminate stabilization term in the gradient-type bilinear form, and contributes nothing to eliminate the stabilization term in the discretization of $L^2$ inner product, which will be discussed in details in next subsection.}}

\textit{Beyond the limitations mentioned above, the gradient projection $Q^K$  satisfying norm equivalence  \eqref{eq:GradProjNormEqui} proposed in the literature is only applicable to eliminating stabilization terms in gradient-type bilinear forms. It does not, however,  address other types of stabilization that may arise in the discrete formulation of the model, such as those associated with the  $L^2$ inner product stemming from a reaction term. These shortcomings motivate the research for a new approach to constructing stabilization-free VEMs.}
We will discuss these issues  in the following subsection.

\subsection{New  strategy in the design of stabilization-free VEMs}
\label{sec::Strategy}

For convenience, we first introduce some necessary notation. Let $\Omega$ be a bounded polytopal domain in $\mathbb R^d\ (d=2,3)$ and $S$ be a subset of $\Omega$. By $\mathbb P_m(S)$ we denote the polynomial space of degree up to $m$ where $m$ is a nonnegative integer. When $m<0$, we set  $\mathbb P_m(S)=\{0\}$. We also use a subspace $\overline{\mathbb P}_m(S)$ consisting of all $m$-degree homogenous polynomials. The symbols $(\cdot,\cdot)_S$ and $\|\cdot\|_S$ stand for the $L^2$ inner product and norm on the square-integrable space  $L^2(S)$, respectively.  $H^1(S)$ is the usual Sobolev space consisting of functions with derivatives in $L^2(S)$ and $H_0^1(S)$ the subspace of functions in $H^1(S)$ with the vanishing trace. $\|\cdot\|_{1,S}$ and $|\cdot|_{1,S}$ are the norm and seminorm on space $H^1(S)$, respectively. When $S=\Omega$, we omit the subscript.

In order to explain the difficulty and strategy in designing stabilization-free VEMs, as an example we consider a simple second-order elliptic problem with reaction term
\begin{equation}\label{eq:2ndEquation}
-\alpha\Delta u+\beta u=f\quad\mbox{in }\Omega,\quad u=0\quad\mbox{on }\partial\Omega,
\end{equation}
where $f\in L^2(\Omega)$, $\alpha$ and $\beta$ are the coefficients of diffusion and reaction terms, respectively. For simplification, we assume that $\alpha$ and $\beta$ are positive constants, bounded above and below.
The corresponding weak form is to find $u\in H_0^1(\Omega)$ such that
\begin{equation*}
\alpha(\nabla u,\nabla v)+\beta(u,v)=( f,v),\quad\forall v\in H_0^1(\Omega).
\end{equation*}

Let $V_h$ be the VE space on a polygonal/polyhedral mesh $\mathcal T_h$. The corresponding VEM is to find $u_h\in V_h$ such that
\begin{equation*}
a_h( u_h, v_h)+b_h(u_h,v_h)=\langle f_h,v_h\rangle,\quad\forall v_h\in V_h.
\end{equation*}
In the common VEMs, the two discrete bilinear forms $a_h( u_h, v_h)$ and $b_h(u_h,v_h)$ should contain the stabilization terms with different factors depending on mesh size, such that they satisfy the stability
\begin{equation}\label{eq:Stability}
a_h^K( v_h, v_h)\sim (\alpha\nabla v_h,\nabla v_h)_K,\quad b_h^K( v_h, v_h)\sim ( \beta v_h,v_h)_K,\quad\forall v_h\in V_h,
\end{equation} 
where $a_h^K( v_h, v_h)$ and $b_h^K( v_h, v_h)$ are the restrictions of $a_h( v_h, v_h)$ and $b_h(v_h,v_h)$ on element $K$ in $\mathcal T_h$, respectively. Through this paper, we use the symbol ``$C$" to denote a general positive constant independent of the mesh size, which may be different at different places. The symbol ``$a\lesssim b$" stands for ``$a\leq Cb$" and ``$a\sim b$" for ``$a\lesssim b\lesssim a$".

The norm equivalence \eqref{eq:GradProjNormEqui} implies that we can use the gradient projection $Q^K$ to design a stabilization-free discrete bilinear form for the gradient term by defining
$$a_h(u_h,v_h)=\sum_{K\in \mathcal T_h}(\alpha Q^K\nabla u_h,Q^K\nabla v_h)_K.$$
\textit{A key question arises:how can one design a stabilization-free bilinear form $b_h(u_h,v_h)$ for the $L^2$ term, in other words, how can one design a stabilization-free scheme that meets the two stability requirements in  \eqref{eq:Stability} simultaneously?
 Obviously, the gradient projection $Q^K$ can not resolve this issue, since it only preserves the norm equivalence \eqref{eq:GradProjNormEqui}.}

We now present a general strategy to eliminate the stabilization terms while  maintaining the stability condition \eqref{eq:Stability} for both resulting discrete bilinear forms  $a_h(\cdot,\cdot)$ and $b_h(\cdot,\cdot)$. The key idea is to construct an interpolation operator  $\widetilde I^K$ from the local VE space into a macro FE space (such as Lagrange FE) defined on a simplicial subtriangulation $\mathcal T_K$ of element $K$, , which satisfies  {\bf norm equivalence}-meaning $\widetilde I^K$ is bounded both above and below with respect to the relevant norms on the VE space
\begin{equation*}
\|v_h\|_K\sim \|\widetilde I^Kv_h\|_K,\quad |v_h|_{1,K}\sim |\widetilde I^Kv_h|_{1,K},\quad v_h\in V_h.
\end{equation*}
This naturally leads us a guide to design the stabilization-free VEM for the general second-order problem with low-order term by simply setting
\begin{equation}\label{eq:ahbhDefinition}
a_h(u_h,v_h)=\sum_{K\in \mathcal T_h}(\alpha\nabla\widetilde I^K u_h,\nabla\widetilde I^K v_h)_K,\quad b_h(u_h,v_h)=\sum_{K\in \mathcal T_h}(\beta\widetilde I^K u_h,\widetilde I^K v_h)_K,
\end{equation}
which satisfy the stability condition \eqref{eq:Stability}.

In addition, the interpolation operator $\widetilde I^K$ must be computable solely from the degrees of freedom of the virtual element space  ({\bf computability}), which ensures the resulting discrete bilinear forms are fully computable. In practice,  by using of such an  operator, we only need to replace the VE function with its interpolation to guarantee the computability of all terms in the discrete scheme.
Furthermore, the interpolation operator  $\widetilde I^K$  
should satisfy two fundamental properties:
\begin{itemize}
\item {\bf Polynomial preservation}:  it preserves all polynomials of degree $\leq k$, where $k$ is the order of the VEM.
\item {\bf Polynomial consistency} : for all  $v_h\in V_h$ and every $q\in\mathbb P_k(K)$,
\begin{equation*}
(\nabla\widetilde I^Kv_h, \nabla q)_K=(\nabla v_h, \nabla q)_K,\quad (\widetilde I^Kv_h,q)_K=(v_h,q)_K.
\end{equation*}
\end{itemize}
% (), and satisfy the {\bf consistency} with respect to polynomials of degree $\leq k$
These conditions ensure that the resulting discrete forms satisfy the standard $k$-consistency, i.e., for any polynomial $q$ of degree $\leq k$, it holds
\begin{equation}\label{eq:Consistency}
a_h^K( v_h, q)=(\alpha\nabla v_h,\nabla q)_K,\quad b_h^K(v_h,q)=(\beta v_h,q)_K,\quad \forall v_h\in V_h.
\end{equation}
As established in \cite{BeiraodaVeiga2013},  the combination of the stability \eqref{eq:Stability} and consistency \eqref{eq:Consistency} guarantees  that the resulting stabilization-free VEM with a well-chosen right hand side $f_h$ admits a unique solution and achieves the expected convergence rate.

Therefore, in designing stabilization-free VEMs based on finite element interpolation, it seems that the interpolation $\widetilde I^K$ should be {\bf computable, polynomial-preserving, norm-equivalent} and {\bf polynomial-consistency}. Following this guideline, we propose a first type of stabilization-free VEM, where the right‑hand side is  taken as the usual  piecewise $L^2$ projection. Because this method retains both the stability \eqref{eq:Stability} and the consistency \eqref{eq:Consistency}, its error analysis follows the same framework as that of the standard VEM.

Through a thorough analysis, we find that the consistency \eqref{eq:Consistency} is in fact not essential for designing stabilization-free VEMs. Motivated by this observation, we introduce a second type of stabilization-free VEM, in which the right hand side is also constructed via finite element interpolation (see Section~\ref{sec:SecondSFVEM}). Notably, this second scheme does not require the underlying FE interpolation to satisfy polynomial consistency; consequently, it does not preserve \eqref{eq:Consistency}.  Its advantages lie in a simpler construction with fewer degrees of freedom, leading to lower computational cost in practice. Furthermore, the second stabilization-free VEM is preferable in solving more complex problems, such as those involving nonlinearities or variable coefficients, where the consistency  \eqref{eq:Consistency} does not hold in general, even if the FE interpolation itself satisfies the polynomial consistency. It is also worth mentioning that the second scheme  gives rise to a new discretization on polytopal meshes, termed the {\bf the projection-based FEM}, see Remark \ref{rmk:PbFEM} for the details.

\textit{We stress that the FE-interpolation-based strategy proposed here can simultaneously remove stabilization terms arising from the discretizations of both diffusion and reaction terms—a feature not achievable by existing gradient‑projection-based techniques. Moreover, the same design philosophy can be extended to other polytopal discretization frameworks, such as the hybrid high-order method~\cite{DiPietro2016} and the weak Galerkin method~\cite{Wang2014,Wang2017}.}

\subsection{Summary of the present work}

This paper proposes a new methodology for eliminating stabilization terms in VEMs by using FE interpolations. Section \ref{sec:SFVEM} introduces an abstract framework that delivers two families of stabilization-free schemes:

\begin{itemize}
\item The first relies on FE interpolations that are computable, polynomial-preserving, norm-equivalent, and polynomial-consistent.
\item The second only demands computability, polynomial-preservation, and norm-equivalence, thereby relaxing the consistency requirement.
\end{itemize}
Abstract error estimates are provided for each family. Then the FE interpolations used in the framework are subsequently constructed for conforming and nonconforming VEMs in 2D and 3D: Section \ref{sec:FirstInterpolation} designs polynomial-consistent interpolations (for first scheme), while Section \ref{sec:SecondInterpolation} builds interpolations without consistency (for second scheme). The latter leads to simpler and more economical constructions. Finally, in Section \ref{sec:NumerTest} we verify the optimal convergence rates of the proposed stabilization-free VEMs and make some comparisons with the standard VEM and FEM through numerical tests. The approach presented here not only removes stabilizers for both diffusion and reaction terms but also extends naturally to other polytopal discretizations, such as the hybrid high-order method and the weak Galerkin method.

\section{The two types of stabilization-free VEMs}
\label{sec:SFVEM}

\subsection{Meshes and notation}

Over the domain $\Omega$, let $\mathcal T_h$ be a polytopal mesh consisting of polygons for $d=2$ and polyhedrons for $d=3$. For a polytope $K\in\mathcal T_h$, $h_K$ is its diameter and $\bm n_K$ is the unit outward normal vector on $\partial K$. The set of vertices of $K$ is denoted by $\mathcal V(K)$, the set of edges of $K$ by $\mathcal E(K)$, and the set of all $(d-1)$-faces of $K$ by $\mathcal F(K)$. Specially for $d=2$, we have $\mathcal E(K)=\mathcal F(K)$.

We assume that there exists a positive constant $\rho$, such that every $K\in\mathcal T_h$ is star-shaped with respect to a disc/ball with radius $\geq\rho h_K$ and the length $|e|$ of every edge $e$ of $K$ is greater than or equal to $\rho h_K$. When $d=3$, each face in $\mathcal T_h$ is star-shaped with respect to the regularity parameter $\rho$ in the same sense. Under this mesh assumption, for each polytope $K$ there exists a shape-regular simplicial subdivision $\mathcal T_K$ with the comparable mesh size relative to $h_K$. Specially for $d=3$, the subdivision $\mathcal T_K$ generates naturally a triangulation $\mathcal T_F$ on face $F$ of $K$. In 2D, we can obtain the simplicial subdivision $\mathcal T_K$ by connecting the center of the disc and vertices of $K$. In 3D, we can also obtain the simplicial subdivision $\mathcal T_K$ by connecting the center of the ball, the center of the disc in each face and vertices of $K$.

For convenience of presentation of our work, we assume that the edges of $K$ are still the edges in the subdivision $\mathcal T_K$, i.e. the edges of $K$ are not subdivided into smaller edges, This is easy to do in practical computations. For the simplex $T\in\mathcal T_K$, the symbols $\mathcal V(T),\mathcal E(T)$ and $\mathcal F(T)$ stand for the sets of vertices, edges and $(d-1)$-faces of $T$, separately. We set 
$$\mathcal V(\mathcal T_K)=\underset{T\in\mathcal T_K}{\cup}\mathcal V(T),\quad\mathcal E(\mathcal T_K)=\underset{T\in\mathcal T_K}{\cup}\mathcal E(T),\quad\mathcal F(\mathcal T_K)=\underset{T\in\mathcal T_K}{\cup}\mathcal F(T).$$
We use $\mathcal V^\partial(\mathcal T_K)$ to denote the set of vertices in $\mathcal V(\mathcal T_K)$ located at boundary of $K$ and $\mathcal V^{\mathrm{int}}(\mathcal T_K)$ the set of vertices in $\mathcal V(\mathcal T_K)$ located at the interior of $K$. Similarly, $\mathcal E^\partial(\mathcal T_K)$ and $\mathcal F^\partial(\mathcal T_K)$ stand for the set of boundary edges and $(d-1)$-faces in $\mathcal T_K$ located on $\partial K$, respectively. $\mathcal E^{\mathrm{int}}(\mathcal T_K)$ and $\mathcal F^{\mathrm{int}}(\mathcal T_K)$ stand for the set of internal edges and $(d-1)$-faces in $\mathcal T_K$ located in the interior of  $K$, respectively. Thus we have the relations
$$\mathcal V(\mathcal T_K)=\mathcal V^\partial(\mathcal T_K)\cup\mathcal V^{\mathrm{int}}(\mathcal T_K),\quad \mathcal E(\mathcal T_K)=\mathcal E^\partial(\mathcal T_K)\cup\mathcal E^{\mathrm{int}}(\mathcal T_K),\quad \mathcal F(\mathcal T_K)=\mathcal F^\partial(\mathcal T_K)\cup\mathcal F^{\mathrm{int}}(\mathcal T_K). $$
Specially for the 2D case ($d=2$), we have $\mathcal V^\partial(\mathcal T_K)=\mathcal V(K)$ and $\mathcal E^\partial(\mathcal T_K)=\mathcal E(K)$ according to the assumption on the subdivision $\mathcal T_K$ of $K$.

For one geometric object $S$ (e.g. an edge, face or element in mesh), $|S|$ stands for its measure, $\bm x_S$ its barycenter and $h_S$ its diameter. For integer $k\geq0$, the set $\mathbb M_k(S)$  of scaled monomials is given by
$$\mathbb M_k(S):=\Big\{\Big(\frac{\bm x-\bm x_S}{h_S}\Big)^{\bm\beta};~|\bm\beta|\leq k\Big\},$$
where $|\bm\beta|=\beta_1+\cdots+\beta_d$, $\bm x^{\bm\beta}=x_1^{\beta_1}\cdots x_d^{\beta_d}$,  $\bm\beta=(\beta_1,\ldots,\beta_d)$ is a nonnegative multi-index, and $d$ is the dimension of $S$. The set $\mathbb M_k(S)$ forms a basis of $k$-degree polynomial space $\mathbb P_k(S)$. 

Let $Q_k^S$ be the $L^2$ projection operator onto polynomial space $\mathbb P_k(S)$.
 Given $v\in H^1(S)$, we define its $H^1$ projection onto $\mathbb P_k(S)$ by finding $\Pi_k^S v\in\mathbb P_k(S)$ such that
\begin{equation*}
\left\{\begin{array}{l}
(\nabla \Pi_k^S v,\nabla q)_S=(\nabla v,\nabla q)_S,\quad\forall q\in\mathbb P_k(S),\\[0.2cm]
Q^{\partial S}_0\Pi_k^S v=Q^{\partial S}_0 v.
\end{array}\right.
\end{equation*}

\subsection{The stabilization-free VEM with consistency}
\label{sec:FirstSFVEM}

Let $V_h$ be the VE space defined on $\mathcal T_h$ with the interpolation operator $I_h$ and $V_h(K)\subseteq H^1(K)$ its restriction to a polytope $K\in\mathcal T_h$. We define the broken $H^1$ norm as 
$$\|v_h\|_{1,h}=\Big(\sum_{K\in\mathcal T_h}\|v_h\|_{1,K}^2\Big)^{1/2},\quad v_h\in V_h.$$
For the conforming and nonconforming VEs, the following inverse inequalities \cite{AyusodeDios2016,BeiraodaVeiga2017a,Brenner2017,Brenner2018,Chen2018,Cangiani2017} hold true
\begin{equation}\label{VEInverseInequality}
\|\Delta v_h\|_K\lesssim h_K^{-1}|v_h|_{1,K},\quad |v_h|_{1,K}\lesssim h_K^{-1}\|v_h\|_K,\quad\forall v_h\in V_h,\ K\in\mathcal T_h.
\end{equation}
Since $Q_0^{\partial K}(v_h-\Pi_k^Kv_h)=0$, the Poincar\'e-Friedrichs inequality, inverse inequality \eqref{VEInverseInequality} and the $H^1$ orthogonality of $\Pi_k^K$ imply the $L^2$ boundedness of $\Pi_k^K$
\begin{equation*}
\|v_h-\Pi_k^Kv_h\|_K\lesssim h_K|v_h-\Pi_k^Kv_h|_{1,K} \lesssim h_K|v_h|_{1,K}\lesssim\|v_h\|_K.
\end{equation*}
Thus, we get the boundedness of $\Pi^K_k$
\begin{equation}\label{PiBound}
|\Pi^K_kv_h|_{m,K}\lesssim |v_h|_{m,K},\quad \forall v_h\in V_h,\quad m=0,1.
\end{equation}

Let $\widetilde V_h(\mathcal T_K)\subseteq H^1(K)$ be a macro FE space defined on the simplicial subdivision $\mathcal T_K$ of a polytope $K$ and $\widetilde I_k^K$ a FE interpolation from $V_h(K)$ to $\widetilde V_h(\mathcal T_K)$. We assume that the interpolation $\widetilde I_k^K$ satisfies the following conditions:
\begin{itemize}
\item[\bf A1] $\widetilde I_k^K$ is computable by the DoFs of VE,
\item[\bf A2]  $\widetilde I_k^K$ is polynomial-preserving
 \begin{equation}
 \widetilde I_k^Kq=q,\ \forall q\in\mathbb P_k(K),
 \end{equation}
 \item[\bf A3]  $\widetilde I_k^K$ satisfies the norm equivalence
\begin{equation}\label{eq:IhNormEqui}
|v_h|_{m,K}\sim|\widetilde I_k^Kv_h|_{m,K},\quad \forall v_h\in V_h(K),\ m=0,1,
\end{equation}

\item[\bf A4]  $\widetilde I_k^K$ is polynomial-consistent: for any given $v_h\in V_h(K)$, it holds
\begin{equation}\label{eq:IhPolyOrth}
(\nabla \widetilde I_k^Kv_h,\nabla q)_K=(\nabla v_h,\nabla q)_K,\quad (\widetilde I_k^Kv_h,q)_K=(v_h,q)_K,\quad \forall q\in\mathbb P_k(K).
\end{equation}
\end{itemize}

Let $\widetilde I_h|_K=\widetilde I_k^K$ for $K\in\mathcal T_h$. With $\widetilde I_k^K$, we define the discrete bilinear forms as in \eqref{eq:ahbhDefinition}. Recalling the discussion in Section \ref{sec::Strategy}, the assumptions {\bf A1-A4} ensure that  the discrete bilinear forms $a_h(\cdot,\cdot)$ and $b_h(\cdot,\cdot)$ satisfy the stability \eqref{eq:Stability} and consistency \eqref{eq:Consistency}. For the right-hand side, one of the common choices is to use the piecewise $L^2$ projection $Q_h$ to define it, where $Q_h|_K:=Q_k^K$ for $K\in\mathcal T_h$. 

We give the first stabilization-free VEM with consistency for the second-order problem \eqref{eq:2ndEquation}, which reads: finding $u_h\in V_h$ such that
\begin{equation}\label{eq:SFVEM1}
a_h( u_h, v_h)+b_h(u_h,v_h)=( f,Q_hv_h),\quad\forall v_h\in V_h.
\end{equation}

By the standard argument like the one used for the usual VEM in \cite{BeiraodaVeiga2013}, we can prove the well-posedness and obtain the error estimates for the stabilization-free VEM with consistency. Here we omit the proof.
\begin{theorem}
Let $u$ be the exact solution to problem \eqref{eq:2ndEquation}, $u_h$ the VE solution to the discrete problem \eqref{eq:SFVEM1}, $u_\pi$ any piecewise $k$-degree polynomial approximation to $u$ on $\mathcal T_h$. It holds
\begin{equation*}
\|u-u_h\|_{1,h}\lesssim \|u-I_hu\|_{1,h}+\|u-u_\pi\|_{1,h}+h\|f-Q_hf\|+E_h,
\end{equation*}
where $E_h$ is the nonconforming error defined by
$$E_h:=\sup_{v_h\in V_h\backslash\{0\}}\left(\frac{1}{\|v_h\|_{1,h}}\sum_{K\in\mathcal T_h}\int_{\partial K}\alpha\frac{\partial u}{\partial\bm n_K} v_h\mathrm ds\right),$$
and the symbol $I_h$ refers to the standard VEM interpolation operator based on the DoFs.
\end{theorem}

\begin{remark}
The nonconforming error $E_h$ vanishes for the conforming VEM. For the nonconforming VEM, the nonconforming error $E_h$ can be estimated with the optimal rate $\mathcal O(h^k)$ due to the weak continuity of VE, refer to \cite{AyusodeDios2016}.
\end{remark}

\subsection{The stabilization-free VEM without consistency}
\label{sec:SecondSFVEM}

The $k$-consistency \eqref{eq:Consistency} is ensured by the polynomial-preservation and polynomial-consistency of FE interpolation.
In fact, the $k$-consistency \eqref{eq:Consistency} is unnecessary in the design and analysis of stabilization-free VEMs based on the FE interpolation, such that we can use the FE interpolation without polynomial consistency, that can make the construction of FE interpolation easier.  In this case, the right-hand side needs to be defined by using the FE interpolation. 

Let $\widetilde W_h(\mathcal T_K)\subseteq H^1(K)$ be another macro FE space defined on the simplicial subdivision of a polytope $K$ and $\widetilde J_k^K$ the corresponding FE interpolation from $V_h(K)$ to $\widetilde W_h(\mathcal T_K)$. We assume that the interpolation $\widetilde J_k^K$ satisfies the following conditions:
\begin{itemize}
\item[\bf B1] $\widetilde J_k^K$ is computable by the DoFs of VE,
\item[\bf B2]  $\widetilde J_k^K$ is polynomial-preserving
 \begin{equation}\label{eq:JPolyPreservation}
 \widetilde J_k^Kq=q,\ \forall q\in\mathbb P_k(K),
 \end{equation}
 \item[\bf B3]  $\widetilde J_k^K$ satisfies the norm equivalence
\begin{equation}\label{eq:JNormEqui}
|v_h|_{m,K}\sim|\widetilde J_k^Kv_h|_{m,K},\quad \forall v_h\in V_h(K),\ m=0,1.
\end{equation}
\end{itemize}
Let $\widetilde J_h|_K:=\widetilde J_k^K$ for $K\in\mathcal T_h$. Then the second stabilization-free VEM without consistency for the second-order problem \eqref{eq:2ndEquation} is designed as: finding $u_h\in V_h$ such that
\begin{equation}\label{eq:SFVEM2}
(\alpha\nabla_h\widetilde J_h u_h, \nabla_h\widetilde J_h v_h)+(\beta\widetilde J_h u_h,\widetilde J_hv_h)=( f,\widetilde J_hv_h),\quad\forall v_h\in V_h,
\end{equation}
where $(\nabla_hv_h)|_K=\nabla(v_h|_K)$ on $K\in\mathcal T_h$.

The norm equivalence \eqref{eq:JNormEqui} ensures immediately the well-posedness of the discrete problem \eqref{eq:SFVEM2}. Further we can obtain the error estimates for the stabilization-free VEM without consistency. 
\begin{theorem}\label{thm:SFVEM2Esti}
Let $u$ be the exact solution to problem \eqref{eq:2ndEquation}, $u_h$ the VE solution to the discrete problem \eqref{eq:SFVEM2}, $u_\pi$ any piecewise $k$-degree polynomial approximation to $u$ on $\mathcal T_h$. It holds
\begin{equation}\label{eq:SFVEM2Esti}
\|u-u_h\|_{1,h}\lesssim \|u-I_hu\|_{1,h}+\|u-u_\pi\|_{1,h}+G_h,
\end{equation}
where $G_h$ is the nonconforming error defined by
$$G_h:=\sup_{v_h\in V_h\backslash\{0\}}\left(\frac{1}{\|v_h\|_{1,h}}\sum_{K\in\mathcal T_h}\int_{\partial K}\alpha\frac{\partial u}{\partial\bm n_K} \widetilde J_hv_h\mathrm ds\right).$$
\end{theorem}

\begin{proof}
Let $w_h=I_hu-u_h$. The norm equivalence \eqref{eq:JNormEqui} implies the coercivity of bilinear forms, so we have
\begin{align*}
C\|w_h\|_{1,h}^2 &\leq (\alpha\nabla_h\widetilde J_h(I_hu-u_h), \nabla_h\widetilde J_h w_h)+(\beta\widetilde J_h (I_hu-u_h),\widetilde J_hw_h)\\
&=(\alpha\nabla_h\widetilde J_hI_hu, \nabla_h\widetilde J_h w_h)+(\beta\widetilde J_h I_hu,\widetilde J_hw_h)-(f,\widetilde J_h w_h).
\end{align*}
Due to the polynomial-preservation property \eqref{eq:JPolyPreservation} of $\widetilde J_h$, we get
\begin{align*}
C\|w_h\|_{1,h}^2 &\leq (\alpha\nabla_h\widetilde J_h(I_hu-u_\pi), \nabla_h\widetilde J_h w_h)+(\beta\widetilde J_h (I_hu-u_\pi),\widetilde J_hw_h)\\
&\quad+(\alpha\nabla_h(u_\pi-u), \nabla_h\widetilde J_h w_h)+(\beta(u_\pi-u),\widetilde J_hw_h) \\
&\quad+(\alpha\nabla u, \nabla_h\widetilde J_h w_h)+(\beta u,\widetilde J_hw_h)-(f,\widetilde J_h w_h) \\
&=(\alpha\nabla_h\widetilde J_h(I_hu-u_\pi), \nabla_h\widetilde J_h w_h)+(\beta\widetilde J_h (I_hu-u_\pi),\widetilde J_hw_h)\\
&\quad+(\alpha\nabla_h(u_\pi-u), \nabla_h\widetilde J_h w_h)+(\beta(u_\pi-u),\widetilde J_hw_h) 
+\sum_{K\in\mathcal T_h}\int_{\partial K}\alpha\frac{\partial u}{\partial\bm n_K}\widetilde J_hw_h\mathrm ds.
\end{align*}
Using the norm equivalence \eqref{eq:JNormEqui}, we get
\begin{equation*}
\|I_hu-u_h\|_{1,h}\lesssim \|u-I_hu\|_{1,h}+\|u-u_\pi\|_{1,h}+\frac{1}{\|w_h\|_{1,h}}\sum_{K\in\mathcal T_h}\int_{\partial K}\alpha\frac{\partial u}{\partial\bm n_K}\widetilde J_hw_h\mathrm ds,
\end{equation*}
which, together with the triangle inequality, leads to the estimate \eqref{eq:SFVEM2Esti}.
\end{proof}

\begin{remark}
We observe there is no the approximation error of the right-hand side for the second stabilization-free VEM \eqref{eq:SFVEM2} based on $\widetilde J_h$. In general, the nonconforming error $G_h$ vanishes for the conforming VEM, if $\widetilde J_hv_h$ is also continuous for the continuous VE function $v_h$. For the nonconforming VEM, the nonconforming error $G_h$ can be estimated with the optimal rate $\mathcal O(h^k)$, if $\widetilde J_hv_h$ is also weakly continuous for the weakly continuous VE function $v_h$. That's indeed the case, i.e., the continuity requirement is naturally met in the construction of interpolation operators. 
\end{remark}

\begin{remark}
By contrast, the first stabilization-free VEM \eqref{eq:SFVEM1} based on $\widetilde I_h$  seems to be simple in the error analysis, because $\widetilde I_h$ satisfies the conditions {\bf A1-A4} such that the resulting discrete bilinear forms satisfy the stability and consistency, which enable the error analysis like the standard VEM.  However,
comparing the constructions of $\widetilde I_h$ and $\widetilde J_h$ in following sections, we observe that the interpolation $\widetilde J_h$ satisfying conditions {\bf B1-B3} is easier to construct with less DoFs and cheap to implement in practical computations.
\end{remark}

\begin{remark}
The FE interpolation $\widetilde I_h$ used in the first stabilization-free VEM \eqref{eq:SFVEM1} need to satisfy the polynomial consistence ({\bf A4}), which makes its construction and implementation somewhat more complex (i.e. add more face or volume bubbles). The two types of stabilization-free VEMs are applicable to more complex problems with nonlinear terms or variable coefficients, but the consistency (1.5) does not hold in general. Therefore, we tend to use the second stabilization-free VEM \eqref{eq:SFVEM2} based on $\widetilde J_h$ due to its simplicity. 
\end{remark}

\begin{remark}
By the usual duality argument, we can also show that the $L^2$ error $\|u-u_h\|$ for the two stabilization-free VEM schemes \eqref{eq:SFVEM1} and \eqref{eq:SFVEM2} is $\mathcal O(h^{k+1})$.
\end{remark}

\section{Interpolation with polynomial consistency}
\label{sec:FirstInterpolation}

\subsection{Interpolation from conforming VE to FE in 2D}  
\setcounter{equation}{0}

In this subsection, we present the construction of 2D interpolation operator $\widetilde I_k^K$ satisfying the conditions {\bf A1-A4}.

\subsubsection{Conforming VE}
Given a polygon $K$, we use the $H^1$ projection $\Pi^K_k$ to define the enhanced conforming VE \cite{Ahmad2013} by setting the local space
\begin{align}\label{2DCVEspace}
V_k^{\mathrm c}(K) =&\Big\{v_h\in H^1(K);\ \Delta v_h\in\mathbb P_k(K),\ v_h|_{\partial K}\in\mathbb B_k(\partial K), \nonumber\\
&\quad (v_h,q)_K=(\Pi_k^K v_h,q)_K, \ q\in \overline{\mathbb P}_k(K)\oplus\overline{\mathbb P}_{k-1}(K)\Big\}
\end{align}
with the DoFs 
\begin{align}
&v(\delta),\quad \delta\in\mathcal V(K), \label{2DCVEDoF1}\\
&\frac{1}{|e|}(v,q)_e,\quad q\in\mathbb M_{k-2}(e),\ e\in\mathcal E(K), \label{2DCVEDoF2}\\
&\frac{1}{|K|}(v,q)_K,\quad q\in\mathbb M_{k-2}(K),\label{2DCVEDoF3}
\end{align}
where the boundary space $\mathbb B_k(\partial K)$ is defined by
$\mathbb B_k(\partial K)=\{v_h\in C^0(\partial K);\ v_h|_e\in\mathbb P_k(e),\ e\in\mathcal E(K)\}$.
Note that the projections $\Pi^K_k$ and $Q^K_k$ are computable by the above DoFs.

\subsubsection{Conforming macro FE}

On a given triangle $T$, we define the conforming FE by setting the local space
$$\widetilde V_k(T)=\mathbb P_k(T)\oplus b_T\overline{\mathbb P}_{k-2}(T)\oplus b_T\overline{\mathbb P}_{k-1}(T)\oplus b_T\overline{\mathbb P}_k(T)$$
with the DoFs
\begin{align}
&v(\delta),\quad\delta\in\mathcal V(T), \label{2DCFEDoF1}\\
& \frac{1}{|e|}(v,q)_e,\quad q\in\mathbb M_{k-2}(e),\ e\in\mathcal E(T), \label{2DCFEDoF2}\\
&\frac{1}{|T|}(v,q)_T,\quad q\in\mathbb M_k(T),\label{2DCFEDoF3}
\end{align}
where $b_T\in\mathbb P_3(T)\cap H_0^1(T)$ is the bubble function obtained by multiplying all the barycentric coordinates on $T$. Following the same discussion as the classical Lagrange FE \cite{Ciarlet1978}, we easily prove the unisolvence of DoFs for this FE. 

On the subtriangulation $\mathcal T_K$ of a polygon $K$, we define the macro FE by setting the space 
$$\widetilde V_k(\mathcal T_K)=\{v_h\in C^0(K);\ v_h|_T\in \widetilde V_k(T),\ T\in\mathcal T_K\}$$
with the above DoFs.

\subsubsection{Interpolation from VE to FE}
\label{sec:Interpolation2D}

For a given VE function $v_h\in V_k^{\mathrm c}(K)$, we define the interpolation $\widetilde I_k^K$ from $V_k^{\mathrm c}(K)$ to $\widetilde V_k(\mathcal T_K)$ by letting $\widetilde I_k^Kv_h\in\widetilde V_k(\mathcal T_K)$ satisfy
\begin{align}
&\widetilde I_k^K v_h(\delta)=v_h(\delta),\quad\delta\in\mathcal V(K), \label{2DCFEIntDoF1}\\
&\widetilde I_k^K v_h(\delta)=\Pi_k^Kv_h(\delta),\quad\delta\in\mathcal V^{\mathrm{int}}(\mathcal T_K), \\
& \frac{1}{| e|}(\widetilde I_k^K v_h-v_h,q)_e=0,\quad q\in\mathbb M_{k-2}(e),\  e\in\mathcal E^\partial(\mathcal T_K),\label{2DCFEIntDoF3}\\
& \frac{1}{|e|}(\widetilde I_k^K v_h-\Pi_k^Kv_h,q)_{ e}=0,\quad q\in\mathbb M_{k-2}( e),\  e\in\mathcal E^{\mathrm{int}}(\mathcal T_K),\\
&\frac{1}{|T|}(\widetilde I_k^Kv_h-Q_k^K v_h,q)_T=0,\quad q\in\mathbb M_k(T),\ T\in\mathcal T_K.\label{2DCFEIntDoF5}
\end{align}

Observing that $v_h$ is a piecewise polynomial of degree $\leq k$ on $\partial K$, $v_h|_{\partial K}$ can be exactly computed by the DoFs \eqref{2DCVEDoF1}-\eqref{2DCVEDoF2}, which, together with the computability of $\Pi_k^K$ and $Q_k^K$, implies the interpolation $\widetilde I_k^K$ is computable by the DoFs \eqref{2DCVEDoF1}-\eqref{2DCVEDoF3} of conforming VE, i.e. it satisfies the condition {\bf A1}.

It can be also seen that $\widetilde I_k^K$ preserves all the polynomials of degree $\leq k$ by recalling the fact that
$$\Pi_k^Kp=Q_k^Kp=p\quad\forall p\in\mathbb P_k(K), \quad \quad \mathbb P_k(K)\subseteq\widetilde V_k(\mathcal T_K) $$
and the unisolvence of DoFs \eqref{2DCFEDoF1}-\eqref{2DCFEDoF3} on space $\widetilde V_k(\mathcal T_K)$. So $\widetilde I_k^K$ satisfies the condition {\bf A2}.

By the interpolation theory of FEM, the trace inequality and inverse inequality, we can obtain the boundedness of $\widetilde I_k^K$.
\begin{lemma}\label{lem:CFEIUpperBound}
For any given $v_h \in V_k^{\mathrm c}(K)$, it holds the boundedness 
\begin{equation}\label{eq:CFEIUpperBound}
|\widetilde I_k^Kv_h|_{m,K}\lesssim |v_h|_{m,K},\quad m=0,1.
\end{equation}
\end{lemma}

\begin{proof}
We first estimate $\|\widetilde I_k^Kv_h-\Pi_k^Kv_h\|_T$ on any given triangle $T$ in $\mathcal T_K$. We observe that the DoFs of $(\widetilde I_k^Kv_h-\Pi_k^Kv_h)$ vanish at internal vertices and internal edges in $\mathcal T_K$, so the definition of $\widetilde I^K_k$ implies
\begin{align}\label{eq:2DCFEIBound1}
(\widetilde I_k^Kv_h-\Pi_k^Kv_h)|_T&=\sum_{\delta\in\mathcal V(T)\cap\mathcal V(K)}(v_h(\delta)-\Pi_k^Kv_h(\delta))\phi_\delta+\sum_{e\in\mathcal E(T)\cap\mathcal E(K)}\sum_{q\in\mathbb M_{k-2}(e)}\frac{1}{|e|}\int_e(v_h-\Pi_k^Kv_h)q\mathrm ds\phi_{e,q}\nonumber\\
&+\sum_{q\in\mathbb M_k(T)}\frac{1}{|T|}\int_T(Q_k^Kv_h-\Pi_k^Kv_h)q\mathrm dx\phi_q,
\end{align}
where $\phi_\delta,\phi_{e,q},\phi_q$ are bounded basis functions associated to the vertices, edges and the interior of $T$, i.e. $\|\phi\|_T\lesssim h_T$.  We only need to estimate the DoFs of $(\widetilde I_k^Kv_h-\Pi_k^Kv_h)$.  

For the DoF at the vertex $\delta \in\mathcal V(T)\cap\mathcal V(K)$, we use the inverse inequality \eqref{VEInverseInequality} and trace inequality to obtain
\begin{align}\label{eq:2DCFEIBound2}
|v_h(\delta)-\Pi_k^Kv_h(\delta)| &\leq \|v_h-\Pi_k^Kv_h\|_{\infty,\partial K}\lesssim h_K^{-\frac{1}2}\|v_h-\Pi_k^Kv_h\|_{\partial K} \nonumber\\
&\lesssim h_K^{-1}\|v_h-\Pi_k^Kv_h\|_K+|v_h-\Pi_k^Kv_h|_{1,K}\nonumber\\
&\lesssim h_K^{-1}\|v_h-\Pi_k^Kv_h\|_K,
\end{align}
where we also used the fact that $v_h$ is piecewise polynomial on $\partial K$.

For the DoFs associated to $q\in\mathbb M_{k-2}(e)$ on the edge $e\in\mathcal E(T)\cap\mathcal E(K)$, the trace inequality, inverse inequality \eqref{VEInverseInequality} and the boundedness of scaling monomials imply
\begin{equation}\label{eq:2DCFEIBound3}
\frac{1}{|e|}\Big|\int_e(v_h-\Pi_k^Kv_h)q\mathrm ds\Big|
\lesssim \frac{1}{|e|}\|v_h-\Pi_k^Kv_h\|_e\|q\|_e
\lesssim h_K^{-1}\|v_h-\Pi_k^Kv_h\|_K.
\end{equation}
For the internal DoFs associated to $q\in\mathbb M_{k-2}(T)$ in the interior of $T$, the boundedness of scaling monomials and the property of $L^2$ projection $Q_k^K$ imply
\begin{align}\label{eq:2DCFEIBound4}
\frac{1}{|T|}\Big|\int_T(Q_k^Kv_h-\Pi_k^Kv_h)q\mathrm dx\Big| &\lesssim h_K^{-2}\|Q_k^Kv_h-\Pi_k^Kv_h\|_T\|q\|_T \nonumber\\
&\lesssim h_K^{-1}\|Q_k^Kv_h-\Pi_k^Kv_h\|_K\lesssim h_K^{-1}\|v_h-\Pi_k^Kv_h\|_K.
\end{align}
Combining \eqref{eq:2DCFEIBound1}-\eqref{eq:2DCFEIBound4} and the boundedness of FE basis functions, we get
$$\|\widetilde I_k^Kv_h-\Pi_k^Kv_h\|_T\lesssim \|v_h-\Pi_k^Kv_h\|_K,$$
which, together with the boundedness \eqref{PiBound} of $\Pi_k^K$, leads to 
$$\|\widetilde I_k^Kv_h-\Pi_k^Kv_h\|_T\lesssim \|v_h\|_K.$$
Then by the triangle inequality, we obtain
\begin{equation}\label{eq:2DCFEIBound6}
\|\widetilde I_k^Kv_h\|_K\leq \|\widetilde I_k^Kv_h-\Pi_k^Kv_h\|_K+\|\Pi_k^Kv_h\|_K\lesssim \|v_h\|_K.
\end{equation}
This completes the proof of the boundedness \eqref{eq:CFEIUpperBound} for $m=0$. 

For the case $m=1$, since $\widetilde I_k^K$ is polynomial-preserving, we use the inverse inequality, the boundedness \eqref{PiBound} of $\Pi_k^K$, Poincar\'e inequality and the estimate \eqref{eq:2DCFEIBound6} to get
\begin{align}\label{eq:2DCFEIBound7}
|\widetilde I_k^Kv_h|_{1,K}&\leq |\widetilde I_k^K(v_h-\Pi_k^Kv_h)|_{1,K}+|\Pi_k^Kv_h|_{1,K}
\lesssim h_K^{-1}\|\widetilde I_k^K(v_h-\Pi_k^Kv_h)\|_K+|v_h|_{1,K} \nonumber\\
&\lesssim h_K^{-1}\|v_h-\Pi_k^Kv_h\|_K+|v_h|_{1,K} \lesssim|v_h|_{1,K}.
\end{align}
The proof is complete.
\end{proof}

Further, the interpolation $\widetilde I_k^K$ satisfies the consistency on the VE space $V_k^{\mathrm c}(K)$, presented in the following lemma.
\begin{lemma}\label{lem:CVE2DH1Orth}
For any $v_h\in V_k^{\mathrm c}(K)$, we have
$$(\nabla \widetilde I_k^Kv_h,\nabla w_h)_K=(\nabla v_h,\nabla w_h)_K,\quad\forall w_h\in V_k^{\mathrm c}(K).$$
\end{lemma}

\begin{proof}
Since $v_h$ is a piecewise polynomial of degree $\leq k$ on $\partial K$, the equations \eqref{2DCFEIntDoF1} and \eqref{2DCFEIntDoF3} imply $\widetilde I_k^K v_h=v_h$ on $\partial K$, which, together with the fact that $\Delta w_h\in\mathbb P_k(K)$ in $K$, leads to
$$(\nabla v_h,\nabla w_h)_K=-(v_h,\Delta w_h)_K+\int_{\partial K}v_h\frac{\partial w_h}{\partial\bm n_K}\mathrm ds=-(\widetilde I_k^K v_h,\Delta w_h)_K+\int_{\partial K}\widetilde I_k^Kv_h\frac{\partial w_h}{\partial\bm n_K}\mathrm ds=(\nabla \widetilde I_k^K v_h,\nabla w_h)_K,$$
where we also used the property \eqref{2DCFEIntDoF5} of $\widetilde I_k^K$ for the first term in the above equation.
\end{proof}

Since $\mathbb P_k(K)\subset V_k^{\mathrm c}(K)$, the polynomial consistency \eqref{eq:IhPolyOrth} of $\widetilde I_k^K$  follows Lemma \ref{lem:CVE2DH1Orth} and the property \eqref{2DCFEIntDoF5} of $\widetilde I_k^K$. By Lemma \ref{lem:CVE2DH1Orth}, we easily obtain the lower boundedness of $\widetilde I_k^K$ in $H^1$ semi-norm
\begin{equation}\label{eq:InterLowBoundH1}
|v_h|_{1,K}\leq |\widetilde I_k^Kv_h|_{1,K},\quad \forall v_h\in V_k^{\mathrm c}(K),
\end{equation}
which, combining the Poincar\'e-Friedrichs inequality and inverse inequality on polynomials, yields the lower bound of $\widetilde I_k^Kv_h$ in $L^2$ norm
\begin{align}\label{eq:InterLowBoundL2}
\|v_h\|_K &\leq \|v_h-Q_0^{\partial K}v_h\|_K+\|Q_0^{\partial K}v_h\|_K =\|v_h-Q_0^{\partial K}v_h\|_K+\|Q_0^{\partial K}\widetilde I_k^Kv_h\|_K \nonumber\\
& \lesssim h_K|v_h|_{1,K}+\|\widetilde I_k^Kv_h\|_K\lesssim h_K|\widetilde I_k^Kv_h|_{1,K}+\|\widetilde I_k^Kv_h\|_K\lesssim\|\widetilde I_k^Kv_h\|_K.
\end{align}
The norm equivalence \eqref{eq:IhNormEqui} of $\widetilde I_k^K$ is immediate following the upper bounds \eqref{eq:CFEIUpperBound} and lower bounds \eqref{eq:InterLowBoundH1}-\eqref{eq:InterLowBoundL2}. Therefore, the 2D interpolation $\widetilde I_k^K$ satisfies the conditions {\bf A3-A4}.

\subsection{Interpolation from conforming VE to FE in 3D}

In this subsection, we present the construction of 3D interpolation operator $\widetilde I_k^K$ satisfying the conditions {\bf A1-A4}.

\subsubsection{Conforming VE}

We present the construction of conforming VE \cite{Ahmad2013} on a polyhedron $K$. For each face $F$ of $K$, $V_k^c(F)$ denotes the face VE space defined on $F$ by \eqref{2DCVEspace} with the computable projections $\Pi_k^F$ and $Q_k^F$. We define the boundary VE space on $\partial K$ by
 $$\mathbb B_k(\partial K):=\{v\in H^1(\partial K);\ v|_F\in V_k^c(F),\ \forall F\in\mathcal F(K)\}.$$
Then we define the 3D conforming VE \cite{Ahmad2013} by setting the local VE space still denoted by $V_k^{\mathrm c}(K)$
 \begin{align*}
 V_k^{\mathrm c}(K)= &\Big\{v_h\in H^1(K);\ \Delta v_h\in\mathbb P_k(K),\ v_h|_{\partial K}\in\mathbb B_k(\partial K) \\
 &\quad (v_h,q)_K=(\Pi_k^K v_h,q)_K, \ q\in \overline{\mathbb P}_k(K)\oplus\overline{\mathbb P}_{k-1}(K)\Big\}
 \end{align*}
 with the DoFs
 \begin{align}
&v(\delta),\quad \delta\in\mathcal V(K) \label{3DCVEDoF1},\\
&\frac{1}{|e|}(v,q)_e,\quad q\in\mathbb M_{k-2}(e),\ e\in\mathcal E(K), \\
&\frac{1}{|F|}(v,q)_F,\quad q\in\mathbb M_{k-2}(F),\ F\in\mathcal F(K), \\
&\frac{1}{|K|}(v,q)_K,\quad q\in\mathbb M_{k-2}(K)\label{3DCVEDoF4}.
\end{align}
Note that the projections $\Pi^K_k$ and $Q^K_k$ are still computable by the above DoFs. 
 
\subsubsection{Conforming macro FE}

On a tetrahedron $T$, we define the conforming FE by setting the local space still denoted by $\widetilde V_k(T)$
$$\widetilde V_k(T)=\mathbb P_k(T)\oplus\sum_{i=0}^3b_T\overline{\mathbb P}_{k-i}(T)\oplus\sum_{F\in\mathcal F(T)}\big(b_F\overline{\mathbb P}_{k-2}(F)\oplus b_F\overline{\mathbb P}_{k-1}(F)\oplus b_F\overline{\mathbb P}_k(F)\big)$$
with the DoFs
\begin{align}
&v(\delta),\quad\delta\in\mathcal V(T), \label{3DCFEDoF1}\\
& \frac{1}{|e|}(v,q)_e,\quad q\in\mathbb M_{k-2}(e),\ e\in\mathcal E(T), \label{3DCFEDoF2}\\
& \frac{1}{|F|}(v,q)_F,\quad q\in\mathbb M_k(F),\ F\in\mathcal F(T), \label{3DCFEDoF3}\\
&\frac{1}{|T|}(v,q)_T,\quad q\in\mathbb M_k(T),\label{3DCFEDoF4}
\end{align}
where $b_T\in\mathbb P_4(T)\cap H_0^1(T)$ is the element bubble function obtained by multiplying all the barycentric coordinates on $T$ and $b_F\in\mathbb P_3(T)$ is the face bubble function obtained by multiplying the three barycentric coordinates with respect to vertices of face $F$. We note that the polynomials in space $b_F\overline{\mathbb P}_m(F)$ can be naturally extended to the tetrahedron $T$ by the Bernstein form in the barycentric coordinate, see \cite{Chen2024a} for the details.
Following the same discussion as the classical Lagrange FE \cite{Ciarlet1978}, we easily prove the unisolvence of DoFs for this FE. 

On the subdivision $\mathcal T_K$ of a polyhedron $K$, we define the macro FE by setting the space 
$$\widetilde V_k(\mathcal T_K)=\{v_h\in C^0(K);\ v_h|_T\in \widetilde V_k(T),\ T\in\mathcal T_K\}$$
with the above DoFs.

\subsubsection{Interpolation from VE to FE}
\label{sec:Interpolation3D}

For a given VE function $v_h\in V_k^{\mathrm c}(K)$, we define the interpolation $\widetilde I_k^K$ from $V_k^{\mathrm c}(K)$ to $\widetilde V_k(\mathcal T_K)$ by letting $\widetilde I_k^Kv_h\in\widetilde V_k(\mathcal T_K)$ satisfy
\begin{align}
&\widetilde I_k^K v_h(\delta)=v_h(\delta),\quad\delta\in\mathcal V (K), \label{3DCFEIntDoF1}\\
&\widetilde I_k^K v_h(\delta)=\Pi_{k,h}^{\partial K}v_h(\delta),\quad\delta\in\mathcal V^{\partial} (\mathcal T_K)\backslash\mathcal V (K),\\
&\widetilde I_k^K v_h(\delta)=\Pi_k^Kv_h(\delta),\quad\delta\in\mathcal V^{\mathrm{int}} (\mathcal T_K),\label{3DCFEIntDoF3}\\
& \frac{1}{| e|}(\widetilde I_k^K v_h-v_h,q)_e=0,\quad q\in\mathbb M_{k-2}(e),\  e\in\mathcal E(K),\label{3DCFEIntDoF4}\\
& \frac{1}{| e|}(\widetilde I_k^K v_h-\Pi_{k,h}^{\partial K}v_h,q)_e=0,\quad q\in\mathbb M_{k-2}(e),\  e\in\mathcal E^\partial(\mathcal T_K)\backslash \mathcal E(K),\label{3DCFEIntDoF5}\\
& \frac{1}{| e|}(\widetilde I_k^K v_h-\Pi_k^Kv_h,q)_e=0,\quad q\in\mathbb M_{k-2}(e),\  e\in\mathcal E^{\mathrm{int}}(\mathcal T_K),\label{3DCFEIntDoF6}\\
& \frac{1}{| F|}(\widetilde I_k^K v_h-Q_{k,h}^{\partial K}v_h,q)_{ F}=0,\quad q\in\mathbb M_k(  F),\  F\in\mathcal F^\partial( \mathcal T_K),\label{3DCFEIntDoF7}\\
& \frac{1}{|F|}(\widetilde I_k^K v_h-\Pi_k^Kv_h,q)_{F}=0,\quad q\in\mathbb M_k( F),\  F\in\mathcal F^{\mathrm{int}}( \mathcal T_K),\\
&\frac{1}{|T|}(\widetilde I_k^Kv_h-Q_k^K v_h,q)_T=0,\quad q\in\mathbb M_k(T),\ T\in\mathcal T_K,\label{3DCFEIntDoF9}
\end{align}
where $\Pi_{k,h}^{\partial K}|_F=\Pi_k^F$ and $Q_{k,h}^{\partial K}|_F=Q_k^F$ for $F\in\mathcal F(K)$.

Similar to the discussion in 2D case, the computability and polynomial-preservation of $\Pi_k^K,\Pi_k^F,Q_k^K$ and $Q_k^F$ imply that the interpolation $\widetilde I_k^K$ is computable by the DoFs \eqref{3DCVEDoF1}-\eqref{3DCVEDoF4} of conforming VE and preserves all polynomial of degree $\leq k$. So $\widetilde I_k^K$ satisfies the conditions {\bf A1-A2}. 

Further, by the interpolation theory of FEM, the trace inequality and inverse inequality, we can obtain the boundedness of $\widetilde I_k^K$.
\begin{lemma}\label{lem:3DCFEIUpperBound}
For any given $v_h \in V_k^{\mathrm c}(K)$, it holds the boundedness 
\begin{equation}\label{eq:3DCFEIUpperBound}
|\widetilde I_k^Kv_h|_{m,K}\lesssim|v_h|_{m,K},\quad m=0,1.
\end{equation}
\end{lemma}

\begin{proof}
We first estimate $\|\widetilde I_k^Kv_h-\Pi_k^Kv_h\|_T$ on a tetrahedron $T$ in $\mathcal T_K$.  Similar to the 2D case,  the DoFs of $(\widetilde I_k^Kv_h-\Pi_k^Kv_h)$ vanish at internal vertices, internal edges and internal faces in $\mathcal T_K$, so the definition of $\widetilde I^K_k$ implies 
\begin{align}\label{eq:3DCFEIUpperBound1}
(\widetilde I_k^Kv_h-\Pi_k^Kv_h)|_T&=\sum_{\delta\in\mathcal V(T)\cap\mathcal V(K)}(v_h(\delta)-\Pi_k^Kv_h(\delta))\phi_\delta+\sum_{\delta\in\mathcal V(T)\cap\mathcal V^\partial(\mathcal T_K)\backslash\mathcal V(K)}(\Pi_{k,h}^{\partial K}v_h(\delta)-\Pi_k^Kv_h(\delta))\phi_\delta \nonumber\\
&+\sum_{e\in\mathcal E(T)\cap\mathcal E(K)}\sum_{q\in\mathbb M_{k-2}(e)}\frac{1}{|e|}\int_e(v_h-\Pi_k^Kv_h)q\mathrm ds\phi_{e,q} \nonumber\\
&+\sum_{e\in\mathcal E(T)\cap\mathcal E^\partial (\mathcal T_K)\backslash\mathcal E(K)}\sum_{q\in\mathbb M_{k-2}(e)}\frac{1}{|e|}\int_e(\Pi_{k,h}^{\partial K}v_h-\Pi_k^Kv_h)q\mathrm ds\phi_{e,q}\nonumber\\
&+\sum_{F\in\mathcal F(T)\cap\mathcal F^{\partial}(\mathcal T_K)}\sum_{q\in\mathbb M_k(F)}\frac{1}{|F|}\int_F(Q_{k,h}^{\partial K}v_h-\Pi_k^Kv_h)q\mathrm ds\phi_{F,q} \nonumber\\
&+\sum_{q\in\mathbb M_k(T)}\frac{1}{|T|}\int_T(Q_k^Kv_h-\Pi_k^Kv_h)q\mathrm dx\phi_q,
\end{align}
where $\phi_\delta,\phi_{e,q},\phi_{F,q},\phi_q$ are bounded basis functions associated to the vertices, edges, faces and the interior of $T$, i.e. $\|\phi\|_T\lesssim h_T^{\frac{3}2}$.  Similar to the 2D case, we only need to estimate the DoFs of $(\widetilde I_k^Kv_h-\Pi_k^Kv_h)$.  

For the DoF at the vertex $\delta$ in $\mathcal V(T)\cap\mathcal V(K)$, the same discussion as in the proof of \eqref{eq:2DCFEIBound2} yields
\begin{equation*}
|v_h(\delta)-\Pi_k^Kv_h(\delta)|\lesssim h_K^{-1}\|v_h-\Pi_k^K v_h\|_{\partial K},
\end{equation*}
since the restriction of $v_h$ to $\partial K$ is the 2D VE function in $\mathbb B_k(\partial K)$ such that the inverse inequality \eqref{VEInverseInequality} holds on $\partial K$. Then we use the inverse inequality \eqref{VEInverseInequality} on 3D VE space and the trace inequality to obtain
\begin{align}
|v_h(\delta)-\Pi_k^Kv_h(\delta)| &\lesssim h_K^{-\frac{3}2}\|v_h-\Pi_k^Kv_h\|_K+h_K^{-\frac{1}2}|v_h-\Pi_k^Kv_h|_{1,K} \nonumber\\
&\lesssim h_K^{-\frac{3}2}\|v_h-\Pi_k^Kv_h\|_K,\quad\quad\delta\in \mathcal V(T)\cap\mathcal V(K).
\end{align}

For the DoF at the vertex $\delta$ in $\mathcal V(T)\cap\mathcal V^\partial(\mathcal T_K)\backslash\mathcal V(K)$, we use the trace inequality, the boundedness \eqref{PiBound} of $\Pi_{k,h}^{\partial K}$ and inverse inequality on polynomials to obtain
$$|\Pi_{k,h}^{\partial K}v_h(\delta)-\Pi_k^Kv_h(\delta)|\lesssim h_K^{-1}\|\Pi_{k,h}^{\partial K}v_h-\Pi_k^Kv_h\|_{\partial K}\lesssim h_K^{-1}(\|v_h\|_{\partial K}+\|\Pi_k^Kv_h\|_{\partial K}),$$
which, together with the trace inequality and the inverse inequality \eqref{VEInverseInequality} on 3D VE space, leads to
\begin{equation}
|\Pi_{k,h}^{\partial K}v_h(\delta)-\Pi_k^Kv_h(\delta)|\lesssim h_K^{-\frac{3}2}(\|v_h\|_K+\|\Pi_k^Kv_h\|_{K}),\quad \delta\in\mathcal V(T)\cap\mathcal V^\partial(\mathcal T_K)\backslash\mathcal V(K).
\end{equation}

For the DoFs on the edge $e\in\mathcal E(T)\cap\mathcal E(K)$, the trace inequality, boundedness of scaling monomials, inverse inequality \eqref{VEInverseInequality} imply
\begin{align}
\frac{1}{|e|}\Big|\int_e(v_h-\Pi_k^Kv_h)q\mathrm ds\Big| &\lesssim\frac{1}{|e|}\|v_h-\Pi_k^Kv_h\|_e\|q\|_e
\lesssim h_K^{-1}\|v_h-\Pi_k^Kv_h\|_{\partial K} \nonumber\\
&\lesssim h_K^{-\frac{3}2}\|v_h-\Pi_k^Kv_h\|_K,\quad q\in\mathbb M_{k-2}(e),\  e\in\mathcal E(T)\cap\mathcal E(K).
\end{align}
For the DoFs on the edge $e\in\mathcal E(T)\cap\mathcal E^\partial(\mathcal T_K)\backslash\mathcal E(K)$, the trace inequality, boundedness of scaling monomials, boundedness \eqref{PiBound} of $\Pi_{k,h}^{\partial K}$ and inverse inequality \eqref{VEInverseInequality} imply
\begin{align}
&\frac{1}{|e|}\Big|\int_e(\Pi_{k,h}^{\partial K}v_h-\Pi_k^Kv_h)q\mathrm ds\Big| \nonumber\\
&\lesssim\frac{1}{|e|}\|\Pi_{k,h}^{\partial K}v_h-\Pi_k^Kv_h\|_e\|q\|_e\lesssim h_K^{-1}\|\Pi_{k,h}^{\partial K}v_h-\Pi_k^Kv_h\|_{\partial K} 
\lesssim h_K^{-1}(\|v_h\|_{\partial K}+\|\Pi_k^Kv_h\|_{\partial K})\nonumber\\
&\lesssim h_K^{-\frac{3}2}(\|v_h\|_K+\|\Pi_k^Kv_h\|_K),\quad q\in\mathbb M_{k-2}(e),\  e\in\mathcal E(T)\cap\mathcal E^\partial(\mathcal T_K)\backslash\mathcal E(K).
\end{align}
For the DoFs on the face $F\in\mathcal F(T)\cap\mathcal F^{\partial}(\mathcal T_K)$, the trace inequality, boundedness of scaling monomials, $L^2$ boundedness of $Q_{k,h}^{\partial K}$ and inverse inequality \eqref{VEInverseInequality} imply
\begin{align}
&\frac{1}{|F|}\Big|\int_F(Q_{k,h}^{\partial K}v_h-\Pi_k^Kv_h)q\mathrm ds\Big| \nonumber\\
 &\leq \frac{1}{|F|}\|Q_{k,h}^{\partial K}v_h-\Pi_k^Kv_h\|_F\|q\|_F \lesssim h_K^{-1}(\|v_h\|_F+\|\Pi_k^Kv_h\|_F) \nonumber\\
&\lesssim h_K^{-\frac{3}2}(\|v_h\|_K+\|\Pi_k^Kv_h\|_K),\quad\quad q\in\mathbb M_k(F),\  F\in\mathcal F(T)\cap\mathcal F^\partial(\mathcal T_K).
\end{align}
For the DoFs in the interior of $T$, the boundedness of scaling monomials and $L^2$ boundedness of $Q_{k}^{K}$ yield
\begin{align}\label{eq:3DCFEIUpperBound7}
\frac{1}{|T|}\Big|\int_T(Q_k^Kv_h-\Pi_k^Kv_h)q\mathrm dx\Big| 
&\leq\frac{1}{|T|}\|Q_k^Kv_h-\Pi_k^Kv_h\|_T\|q\|_T\lesssim h_K^{-\frac{3}2}\|v_h-\Pi_k^Kv_h\|_T,\quad   q\in\mathbb M_k(T).
\end{align}
Combining \eqref{eq:3DCFEIUpperBound1}-\eqref{eq:3DCFEIUpperBound7} and the boundedness of FE basis functions and $\Pi_k^K$, we get the upper bound \eqref{eq:3DCFEIUpperBound} with $m=0$. For $m=1$, the same discussion as in the proof of \eqref{eq:2DCFEIBound7} leads to the upper bound \eqref{eq:3DCFEIUpperBound}. The proof is complete.
\end{proof}

Different from the 2D case (see Lemma \ref{lem:CVE2DH1Orth}), for the 3D case we only obtain the polynomial consistency of $\widetilde I_k^K$ as presented in the following lemma, since the VE function $v_h$ is not piecewise polynomial on $\partial K$ in 3D.
\begin{lemma}\label{lem:3DInterOrthogonality}
For any given $v_h\in V_k^{\mathrm c}(K)$, it holds
\begin{equation}\label{eq:3DInterOrthogonality}
(\nabla \widetilde I_k^Kv_h,\nabla q)_K=(\nabla v_h,\nabla q)_K,\quad (\widetilde I_k^Kv_h,q)_K=(v_h,q)_K,\quad \forall q\in\mathbb P_k(K).
\end{equation}

\begin{proof}
By integration by parts, we use the properties \eqref{3DCFEIntDoF7} and \eqref{3DCFEIntDoF9} of $\widetilde I_k^K$ to obtain
$$(\nabla \widetilde I_k^Kv_h,\nabla q)_K=-(\widetilde I_k^Kv_h,\Delta q)_K+\int_{\partial K}\widetilde I_k^Kv_h\frac{\partial q}{\partial \bm n_K}\mathrm ds=-(v_h,\Delta q)_K+\int_{\partial K}v_h\frac{\partial q}{\partial \bm n_K}\mathrm ds=(\nabla v_h,\nabla q)_K,$$
which is the first one in \eqref{eq:3DInterOrthogonality}. The second one in \eqref{eq:3DInterOrthogonality} follows from the property \eqref{3DCFEIntDoF9} of $\widetilde I_k^K$.
\end{proof}

\end{lemma}
Next we present the lower boundedness of $\widetilde I_k^K$ in the following lemma.
\begin{lemma}\label{lem:3DinterLowBound}
For $\widetilde I_k^K$, it holds
\begin{equation}\label{eq:3DinterLowBound}
|v_h|_{m,K} \lesssim|\widetilde I_k^Kv_h|_{m,K},\quad\forall v_h\in V_k^{\mathrm c}(K),\quad m=0,1.
\end{equation}
\end{lemma}

\begin{proof}
We observe that $V_k^c(F)$ is the restriction of $V_k^{\mathrm c}(K)$ onto face $F$ of $K$ and $\widetilde V_k(\mathcal T_F)$ is the restriction of $\widetilde V_k(\mathcal T_K)$ on the subtriangulation $\mathcal T_F$ of $F$. Taking a closer look, we have 
$$\widetilde I_k^Kv_h|_{\partial K}=\widetilde I_k^{\partial K}v_h,\quad\forall v_h\in V_k^{\mathrm c}(K),$$
since they have the same DoFs on $\partial K$, where $\widetilde I_k^{\partial K}|_F=\widetilde I_k^F$ is the interpolation from VE space $V_k^c(F)$ to FE space $\widetilde V_k(\mathcal T_F)$ defined by \eqref{2DCFEIntDoF1}-\eqref{2DCFEIntDoF5} on faces of $K$.
For $\widetilde I_k^{\partial K}$, we have proven the norm equivalence \eqref{eq:IhNormEqui}, so we have the lower bounds
\begin{equation}\label{eq:BInterLowBoundH1}
|v_h|_{m,\partial K}\lesssim |\widetilde I_k^{\partial K}v_h|_{m,\partial K},\quad m=0,1.
\end{equation}

With the above preparations, we start to prove \eqref{eq:3DinterLowBound} with $m=1$. First by integration by parts, we obtain
\begin{align}\label{eq:3DGreenFormula}
|v_h|_{1,K}^2 &=|v_h-Q_0^{\partial K}v_h|_{1,K}^2=-(v_h-Q_0^{\partial K}v_h,\Delta v_h)_K+\int_{\partial K}(v_h-Q_0^{\partial K}v_h)\frac{\partial v_h}{\partial\bm n_K}\mathrm ds \nonumber\\
&=-(\widetilde I_k^Kv_h-Q_0^{\partial K}v_h,\Delta v_h)_K+\int_{\partial K}(v_h-Q_0^{\partial K}v_h)\frac{\partial v_h}{\partial\bm n_K}\mathrm ds,
\end{align}
where we used the property \eqref{3DCFEIntDoF9} of $\widetilde I_k^K$ for the first term in the above equation.

Observing $Q_0^{\partial K}\widetilde I_k^Kv_h=Q_0^{\partial K} v_h$, for the first term in the above equation, we use the Poincar\'e-Friedrichs inequality and the inverse inequality \eqref{VEInverseInequality} to obtain
\begin{align}\label{eq:3DGreenFormula1}
(\widetilde I_k^Kv_h-Q_0^{\partial K}v_h,\Delta v_h)_K &=(\widetilde I_k^Kv_h-Q_0^{\partial K}(\widetilde I_k^Kv_h),\Delta v_h)_K\nonumber\\
&\leq\|\widetilde I_k^Kv_h-Q_0^{\partial K}(\widetilde I_k^Kv_h)\|_K\|\Delta v_h\|_K \lesssim |\widetilde I_k^Kv_h|_{1,K}|v_h|_{1,K}.
\end{align}
For the second term in \eqref{eq:3DGreenFormula}, we have
\begin{equation}\label{eq:3DGreenFormula2}
\int_{\partial K}(v_h-Q_0^{\partial K}v_h)\frac{\partial v_h}{\partial\bm n_K}\mathrm ds\leq\|v_h-Q_0^{\partial K}v_h\|_{\frac{1}2,\partial K}\Big\|\frac{\partial v_h}{\partial\bm n_K}\Big\|_{-\frac{1}2,\partial K}.
\end{equation}
Following the discussion in \cite[Section 2.5]{Brenner2018}, we get 
$$\|v_h-Q_0^{\partial K}v_h\|_{\frac{1}2,\partial K}\lesssim h_K^{-\frac{1}2}\|v_h-Q_0^{\partial K}v_h\|_{\partial K}+h_K^{\frac{1}2}|v_h|_{1,\partial K}\lesssim Ch_K^{\frac{1}2}|v_h|_{1,\partial K},$$
which, together with the inequality \eqref{eq:BInterLowBoundH1}, trace inequality and inverse inequality, leads to
$$\|v_h-Q_0^{\partial K}v_h\|_{\frac{1}2,\partial K}\lesssim h_K^{\frac{1}2}|v_h|_{1,\partial K}\lesssim h_K^{\frac{1}2}|\widetilde I_k^{\partial K}v_h|_{1,\partial K}\lesssim |\widetilde I_k^Kv_h|_{1,K},$$
where we also used the fact that $\widetilde I_k^Kv_h=\widetilde I_k^{\partial K}v_h$ on $\partial K$.

For $\big\|\frac{\partial v_h}{\partial\bm n_K}\big\|_{-\frac{1}2,\partial K}$, observing $\nabla v_h$ belongs to the space $\bm H(\mathrm{div};K)$, we use the trace inequality on space $\bm H(\mathrm{div};K)$ (see \cite[inequality (2.13)]{BeiraodaVeiga2022}) and the inverse inequality \eqref{VEInverseInequality} to obtain
$$\Big\|\frac{\partial v_h}{\partial\bm n_K}\Big\|_{-\frac{1}2,\partial K}\lesssim |v_h|_{1,K}+h_K\|\Delta v_h\|_K\lesssim |v_h|_{1,K}.$$
Substituting the two estimates above into inequality \eqref{eq:3DGreenFormula2}, we get
\begin{equation}\label{eq:3DGreenFormula3}
\int_{\partial K}(v_h-Q_0^{\partial K}v_h)\frac{\partial v_h}{\partial\bm n_K}\mathrm ds\lesssim |\widetilde I_k^Kv_h|_{1,K}|v_h|_{1,K}.
\end{equation}
Substituting the estimates \eqref{eq:3DGreenFormula1} and \eqref{eq:3DGreenFormula3} into \eqref{eq:3DGreenFormula} yields the lower bound \eqref{eq:3DinterLowBound} with $m=1$. For $m=0$, the same discussion as in the proof of \eqref{eq:InterLowBoundL2} yields the  lower bound \eqref{eq:3DinterLowBound}.
\end{proof}

According to Lemmas \ref{lem:3DCFEIUpperBound}-\ref{lem:3DinterLowBound}, the 3D interpolation $\widetilde I_k^K$ satisfies the conditions {\bf A3-A4}.

\subsection{Interpolation from nonconforming VE to FE}  

In this subsection, we present the construction of interpolation operator $\widetilde I_k^{\mathrm{nc},K}$ from nonconforming VE to FE satisfying the conditions {\bf A1-A4}.
\subsubsection{Nonconforming VE}

By the utilization of $\Pi_k^K$, we define the local space for the nonconforming VE by, for $d=2$
 \begin{align*}
 V_k^{\mathrm{nc}}(K)=\Big\{v\in V_{k+1}^c(K);~(v,q)_K=(\Pi_k^K v,q)_K,\ \forall q\in\overline{\mathbb P}_{k-1}(K), \ v(\delta)=\Pi_k^K v(\delta),\ \forall\delta\in\mathcal V(K)\Big\},
 \end{align*}
 and for $d=3$
 \begin{align*}
 V_k^{\mathrm{nc}}(K)=\Big\{v\in V_{k+1}^c(K);~(v,q)_K=(\Pi_k^Kv,q)_K,\ \forall q\in\overline{\mathbb P}_{k-1}(K), \ v|_e=(\Pi_k^K v)|_e, \forall e\in\mathcal E(K)\Big\}.
 \end{align*}
 Note that for $v_h\in V_k^{\mathrm{nc}}(K)$ in 3D case, we still have $v_h(\delta)=\Pi_k^Kv_h(\delta)$ at vertex $\delta$ of $K$ from the definition.
The  DoFs are defined by 
\begin{align}
&\frac{1}{|F|}(v,q)_F,\quad q\in\mathbb M_{k-1}(F),\ F\in\mathcal F(K), \label{NVEDoF1}\\
&\frac{1}{|K|}(v,q)_K,\quad q\in\mathbb M_{k-2}(K).\label{NVEDoF2}
\end{align}

We can verify that the $H^1$ projection $\Pi_k^K$ from $V_k^{\mathrm{nc}}(K)$ onto $\mathbb P_k(K)$ are computable only by the above DoFs  of nonconforming VE. In fact, the integration by parts yields
 $$(\nabla v_h,\nabla q)_K=-(v_h,\Delta q)_K+\int_{\partial K}v_h\frac{\partial q}{\partial\bm n_K}\mathrm ds.$$
 It is easy to see that for $v_h\in V_k^{\mathrm{nc}}(K)$ and $q\in\mathbb P_k(K)$, $(\nabla v_h,\nabla q)_K$ can be computed only by the DoFs \eqref{NVEDoF1}-\eqref{NVEDoF2}, since $\Delta q$ is a polynomial of order up to $k-2$ in $K$ and $\frac{\partial q}{\partial\bm n_K}$ is a polynomial of order up to $(k-1)$ on each $(d-1)$-face of $K$. Then, similar arguments in \cite{Ahmad2013} yield the unisolvence of  DoFs for space $V_k^{\mathrm{nc}}(K)$.  From the definition of space $V_k^{\mathrm{nc}}(K)$, it is easy to see that the $L^2$ projection $Q_{k+1}^K$ is also computed by these DoFs. We note that for $v_h\in V_k^{\mathrm{nc}}(K)$ in 3D, both the $H^1$ projection $\Pi_{k+1}^F$ and $L^2$ projection $Q_{k+1}^F$ on face $F$ of $K$ are also computable by the above DoFs, since $v_h=\Pi_k^Kv_h$ on $\partial F$.

\subsubsection{Interpolation from nonconforming VE to FE in 2D}

Let $K$ be a polygon. We use the macro FE space $\widetilde V_{k+1}(\mathcal T_K)$ of degree $(k+1)$ as the interpolation space to define the  interpolation from the nonconforming VE to FE. For a given VE function $v_h\in V_k^{\mathrm{nc}}(K)$, we define the interpolation $\widetilde I_k^{\mathrm{nc},K}$ from $V_k^{\mathrm{nc}}(K)$ to $\widetilde V_{k+1}(\mathcal T_K)$ by letting $\widetilde I_k^{\mathrm{nc},K}v_h\in\widetilde V_{k+1}(\mathcal T_K)$ satisfy
\begin{align}
&\widetilde I_k^{\mathrm{nc},K} v_h(\delta)=v_h(\delta),\quad\delta\in\mathcal V(K), \label{2DJhInterProp1}\\
&\widetilde I_k^{\mathrm{nc},K} v_h(\delta)=\Pi_k^Kv_h(\delta),\quad\delta\in\mathcal V^{\mathrm{int}}(\mathcal T_K), \label{2DJhInterProp2}\\
& \frac{1}{|e|}(\widetilde I_k^{\mathrm{nc},K} v_h-v_h,q)_e=0,\quad q\in\mathbb M_{k-1}(e),\ e\in\mathcal E^\partial(\mathcal T_K), \label{2DJhInterProp3}\\
& \frac{1}{|e|}(\widetilde I_k^{\mathrm{nc},K} v_h-\Pi_k^Kv_h,q)_e=0,\quad q\in\mathbb M_{k-1}(e),\ e\in\mathcal E^{\mathrm{int}}(\mathcal T_K), \\
&\frac{1}{|T|}(\widetilde I_k^{\mathrm{nc},K}v_h-Q_{k+1}^K v_h,q)_T=0,\quad q\in\mathbb M_{k+1}(T),\ T\in\mathcal T_K.
\end{align}
Due to the properties of $\Pi_k^K$ and $Q_{k+1}^K$, we can verify that the interpolation $\widetilde I_k^{\mathrm{nc},K}$ is computable by the DoFs \eqref{NVEDoF1}-\eqref{NVEDoF2} of nonconforming VE and preserves all polynomials of degree $\leq k$ on $K$. 

Similar to the discussion in the proof of Lemma \ref{lem:CFEIUpperBound}, by the interpolation theory of FEM, the trace inequality and inverse inequality, we can obtain the boundedness of $\widetilde I_k^{\mathrm{nc},K}$
\begin{equation*}
|\widetilde I_k^{\mathrm{nc},K}v_h|_{m,K}\lesssim |v_h|_{m,K},\quad v_h \in V_k^{\mathrm{nc}}(K),\quad m=0,1.
\end{equation*}
Further, the interpolation $\widetilde I_k^{\mathrm{nc},K}$ has the consistency on the VE space $V_k^{\mathrm{nc}}(K)$.
\begin{lemma}\label{lem:NCVE2DH1Orth}
For any given $v_h\in V_k^{\mathrm{nc}}(K)$, we have
\begin{equation*}
(\nabla \widetilde I_k^{\mathrm{nc},K}v_h,\nabla w_h)_K=(\nabla v_h,\nabla w_h)_K,\quad \forall w_h\in V_k^{\mathrm{nc}}(K).
\end{equation*}
\end{lemma}

\begin{proof}
By integration by parts, we get
$$(\nabla v_h,\nabla w_h)_K=-(v_h,\Delta w_h)_K+\int_{\partial K}v_h\frac{\partial w_h}{\partial\bm n_K}\mathrm ds.$$
Observing $\Delta w_h\in\mathbb P_{k+1}(K)$ in $K$, we have
$$(v_h,\Delta w_h)_K=(Q_{k+1}^Kv_h,\Delta w_h)_K=\sum_{T\in\mathcal T_K}(Q_{k+1}^Kv_h,\Delta w_h)_T=\sum_{T\in\mathcal T_K}(\widetilde I_k^{\mathrm{nc},K}v_h,\Delta w_h)_T=(\widetilde I_k^{\mathrm{nc},K}v_h,\Delta w_h)_K.$$
Since $v_h$ is a piecewise $(k+1)$-degree polynomial on $\partial K$, the properties \eqref{2DJhInterProp1} and \eqref{2DJhInterProp3} of $\widetilde I_k^{\mathrm{nc},K}$ imply $v_h=\widetilde I_k^{\mathrm{nc},K}v_h$ on $\partial K$, so we have
$$\int_{\partial K}v_h\frac{\partial w_h}{\partial\bm n_K}\mathrm ds=\int_{\partial K}\widetilde I_k^{\mathrm{nc},K}v_h\frac{\partial w_h}{\partial\bm n_K}\mathrm ds.$$
Then we have
$$(\nabla v_h,\nabla w_h)_K=-(\widetilde I_k^{\mathrm{nc},K} v_h,\Delta w_h)_K+\int_{\partial K}\widetilde I_k^{\mathrm{nc},K}v_h \frac{\partial w_h}{\partial\bm n_K}\mathrm ds=(\nabla\widetilde I_k^{\mathrm{nc},K} v_h,\nabla w_h)_K.$$
The proof is complete.
\end{proof}

By Lemma \ref{lem:NCVE2DH1Orth} and the similar argument as \eqref{eq:InterLowBoundL2}, we easily obtain the lower boundedness of $\widetilde I_k^{\mathrm{nc},K}$
\begin{equation*}
|v_h|_{m,K}\leq |\widetilde I_k^{\mathrm{nc},K}v_h|_{m,K},\quad \forall v_h\in V_k^{\mathrm{nc}}(K),\ m=0,1.
\end{equation*}
Collecting the results above, we know that the interpolation $\widetilde I_k^{\mathrm{nc},K}$ satisfies the conditions {\bf A1-A4}.

\subsubsection{Interpolation from nonconforming VE to FE in 3D}

Let $K$ be a polyhedron.
For a given VE function $v_h\in V_k^{\mathrm{nc}}(K)$, we define the interpolation $\widetilde I_k^{\mathrm{nc},K}$ from $V_k^{\mathrm{nc}}(K)$ to $\widetilde V_{k+1}(\mathcal T_K)$ by letting $\widetilde I_k^{\mathrm{nc},K}v_h\in\widetilde V_{k+1}(\mathcal T_K)$ satisfy
\begin{align}
&\widetilde I_k^{\mathrm{nc},K} v_h(\delta)=v_h(\delta),\quad\delta\in\mathcal V(K), \\
&\widetilde I_k^{\mathrm{nc},K} v_h(\delta)=\Pi_{k+1,h}^{\partial K} v_h(\delta),\quad\delta\in\mathcal V^\partial(\mathcal T_K)\backslash\mathcal V(K), \\
&\widetilde I_k^{\mathrm{nc},K} v_h(\delta)=\Pi_k^K v_h(\delta),\quad\delta\in\mathcal V^{\mathrm{int}}(\mathcal T_K), \\
& \frac{1}{|e|}(\widetilde I_k^{\mathrm{nc},K} v_h-v_h,q)_e=0,\quad q\in\mathbb M_{k-1}(e),\ e\in\mathcal E(K), \\
& \frac{1}{|e|}(\widetilde I_k^{\mathrm{nc},K} v_h-\Pi_{k+1,h}^{\partial K}v_h,q)_e=0,\quad q\in\mathbb M_{k-1}(e),\ e\in\mathcal E^\partial(\mathcal T_K)\backslash\mathcal E( K), \\
& \frac{1}{|e|}(\widetilde I_k^{\mathrm{nc},K} v_h-\Pi_k^Kv_h,q)_e=0,\quad q\in\mathbb M_{k-1}(e),\ e\in\mathcal E^{\mathrm{int}}(\mathcal T_K),\\
&\frac{1}{| F|}(\widetilde I_k^{\mathrm{nc},K} v_h-Q_{k+1,h}^{\partial K}v_h,q)_{ F}=0,\quad q\in\mathbb M_{k+1}( F),\  F\in\mathcal F^{\partial}(\mathcal T_K),\label{3DJInterpolationEq7}\\
&\frac{1}{|F|}(\widetilde I_k^{\mathrm{nc},K} v_h-\Pi_k^Kv_h,q)_F=0,\quad q\in\mathbb M_{k+1}(F),\ F\in\mathcal F^{\mathrm{int}}(\mathcal T_K),\\
&\frac{1}{|T|}(\widetilde I_k^{\mathrm{nc},K}v_h-Q_{k+1}^K v_h,q)_T=0,\quad q\in\mathbb M_{k+1}(T),\ T\in\mathcal T_K.\label{3DJInterpolationEq9}
\end{align}

Due to the property of $\Pi_{k+1}^F,\Pi_k^K, Q_{k+1}^F$ and $Q_{k+1}^K$, we can verify that the interpolation $\widetilde I_k^{\mathrm{nc},K}$ is computable by the DoFs \eqref{NVEDoF1}-\eqref{NVEDoF2} of nonconforming VE and preserves all polynomials of degree $\leq k$ on $K$.

By the interpolation theory of FEM, the trace inequality and inverse inequality, we can obtain the upper boundedness of $\widetilde I_k^{\mathrm{nc},K}$, i.e.
\begin{equation}\label{eq:3DJInterUpperBound}
|\widetilde I_k^{\mathrm{nc},K}v_h|_{m,K}\lesssim |v_h|_{m,K},\quad\forall v_h\in V_k^{\mathrm{nc}}(K),\ m=0,1.
\end{equation}

We present the lower boundedness of $|\widetilde I_k^{\mathrm{nc},K}v_h|_{1,K}$.
\begin{lemma}\label{lem:3DJInterLowBound}
For $\widetilde I_k^{\mathrm{nc},K}$, it holds
$$|v_h|_{m,K}\lesssim|\widetilde I_k^{\mathrm{nc},K}v_h|_{m,K},\quad\forall v_h\in V_k^{\mathrm{nc}}(K),\ m=0,1.$$
\end{lemma}

\begin{proof}
We observe that the restriction of $V_k^{\mathrm{nc}}(K)$ onto face $F$ of $K$ is a subspace of the conforming VE space $V_{k+1}^c(F)$ on face $F$ and $\widetilde V_{k+1}(\mathcal T_F)$ is the restriction of $\widetilde V_{k+1}(\mathcal T_K)$ on the subtriangulation $\mathcal T_F$ of $F$. Taking a closer look, we have 
$$\widetilde I_k^{\mathrm{nc},K}v_h|_{\partial K}=\widetilde I_{k+1}^{\partial K}v_h,\quad\forall v_h\in V_k^{\mathrm c}(K),$$
because $\widetilde I_k^{\mathrm{nc},K}v_h$ and $\widetilde I_{k+1}^{\partial K}v_h$ have the same DoFs on $\partial K$, where $\widetilde I_{k+1}^{\partial K}|_F=\widetilde I_{k+1}^F$ is the interpolation from VE space $V_{k+1}^c(F)$ to FE space $\widetilde V_{k+1}(\mathcal T_F)$ defined by \eqref{2DCFEIntDoF1}-\eqref{2DCFEIntDoF5} with $(k+1)$ instead of $k$ on faces of $K$. The remaining proof is completely similar to that of Lemma \ref{lem:3DinterLowBound}, so we omit it.
\end{proof}

Collecting the results above, we know that the interpolation $\widetilde I_k^K$ satisfies the conditions {\bf A1-A4}.

\section{ Interpolation without polynomial consistency}
\label{sec:SecondInterpolation}
\setcounter{equation}{0}

In this section, we shall present these FE interpolations without polynomial consistency but satisfying the conditions {\bf B1-B3} for the conforming and nonconforming VEs in 2D and 3D.

\subsection{Interpolation from conforming VE to FE in 2D}

Let $K$ be a polygon with the subdivision $\mathcal T_K$. Given a triangle $T$ in $\mathcal T_K$, we define the conforming FE by setting the local space
$$\widetilde W_k(T)=\mathbb P_k(T)\oplus b_T\overline{\mathbb P}_{k-2}(T)$$
with the DoFs
\begin{align}
&v(\delta),\quad\delta\in\mathcal V(T), \label{W2DCFEDoF1}\\
& \frac{1}{|e|}(v,q)_e,\quad q\in\mathbb M_{k-2}(e),\ e\in\mathcal E(T), \label{W2DCFEDoF2}\\
&\frac{1}{|T|}(v,q)_T,\quad q\in\mathbb M_{k-2}(T).\label{W2DCFEDoF3}
\end{align}
On the subdivision $\mathcal T_K$, we define the macro FE by setting the space 
$$\widetilde W_k(\mathcal T_K)=\{v_h\in C^0(K);\ v_h|_T\in \widetilde W_k(T),\ T\in\mathcal T_K\}$$
with the above DoFs.

For $v_h\in V_k^{\mathrm c}(K)$, we define the interpolation $\widetilde J_k^K$ from $V_k^{\mathrm c}(K)$ to $\widetilde W_k(\mathcal T_K)$ by letting $\widetilde J_k^Kv_h\in\widetilde W_k(\mathcal T_K)$ satisfy
\begin{align}
&\widetilde J_k^K v_h(\delta)=v_h(\delta),\quad\delta\in\mathcal V (K), \label{W2DCFEIntDoF1}\\
&\widetilde J_k^K v_h(\delta)=\Pi_k^Kv_h(\delta),\quad\delta\in\mathcal V^{\mathrm{int}}(\mathcal T_K), \\
& \frac{1}{| e|}(\widetilde J_k^K v_h-v_h,q)_e=0,\quad q\in\mathbb M_{k-2}(e),\  e\in\mathcal E^\partial(\mathcal T_K),\label{W2DCFEIntDoF3}\\
& \frac{1}{|e|}(\widetilde J_k^K v_h-\Pi_k^Kv_h,q)_{ e}=0,\quad q\in\mathbb M_{k-2}( e),\  e\in\mathcal E^{\mathrm{int}}(\mathcal T_K),\\
&\frac{1}{|T|}(\widetilde J_k^Kv_h-Q_k^K v_h,q)_T=0,\quad q\in\mathbb M_{k-2}(T),\ T\in\mathcal T_K.\label{W2DCFEIntDoF5}
\end{align}
Similar to the previous discussion in Section \ref{sec:Interpolation2D}, we easily verify the computability and polynomial-preservation of $\widetilde J_k^K$, i.e. satisfying the conditions {\bf B1-B2}. By the interpolation theory of FEM, the same discussion as in the proof of Lemma \ref{lem:CFEIUpperBound} yields the upper bound
\begin{equation}\label{widetildeJUB}
|\widetilde J_k^Kv_h|_{m,K}\lesssim|v_h|_{m,K},\quad \forall v_h\in V_k^{\mathrm c}(K),\ m=0,1.
\end{equation}
Before proving the lower bound of $\widetilde J_k^Kv_h$, we observe 
\begin{equation}\label{eq:PiIdentity}
\Pi_k^Kv_h=\Pi_k^K(\widetilde J_k^Kv_h),\quad v_h\in V_k^{\mathrm c}(K).
\end{equation}
This is easily verified by the similar argument in proof of Lemma \ref{lem:CVE2DH1Orth}. Next we show the lower bound of $\widetilde J_k^Kv_h$, which, together with \eqref{widetildeJUB},  implies that $\widetilde J_k^Kv_h$ satisfies the condition {\bf B3}.
\begin{lemma}
For any given $v_h \in V_k^{\mathrm c}(K)$, it holds the boundedness 
\begin{equation}\label{eq:WCFEILowerBound}
|v_h|_{m,K}\lesssim|\widetilde J_k^Kv_h|_{m,K},\quad m=0,1.
\end{equation}
\end{lemma}

\begin{proof}
We let $\#\mathcal V(K)$ denote the number of vertices of $K$ and set
$$P_{K}v_h=\frac{1}{\#\mathcal V(K)}\sum_{\delta\in\mathcal V(K)}v_h(\delta).$$
It is immediate that 
\begin{equation}\label{eq:JLowerBound0}
P_{K}v_h=P_{K}(\widetilde J_k^Kv_h),\quad \|\widetilde J_k^Kv_h-P_{K}(\widetilde J_k^Kv_h)\|_K\lesssim h_K|\widetilde J_k^Kv_h|_{1,K}.
\end{equation}
Then the definition of $\widetilde J_k^K$ implies that
\begin{align}\label{eq:JLowerBound1}
|v_h|_{1,K}^2 &=-(v_h,\Delta v_h)_K+\int_{\partial K}v_h\frac{\partial v_h}{\partial\bm n_K}\mathrm ds \nonumber\\
&=-(v_h,\Delta v_h)_K+\int_{\partial K}\widetilde J_k^Kv_h\frac{\partial v_h}{\partial\bm n_K}\mathrm ds \nonumber\\
&=(\nabla \widetilde J_k^K v_h,\nabla v_h)_K+(\widetilde J_k^K v_h,\Delta v_h)_K-(v_h,\Delta v_h)_K\nonumber\\
&=(\nabla \widetilde J_k^K v_h,\nabla v_h)_K+(\widetilde J_k^K v_h-P_K(\widetilde J_k^K v_h),\Delta v_h)_K-(v_h-P_K(\widetilde J_k^K v_h),\Delta v_h)_K.
\end{align}
For the last term in \eqref{eq:JLowerBound1}, we observe that $\Delta v_h$ is a polynomial of degree $\leq k$, so there exists a unique decomposition
$$\Delta v_h=q_1+q_2, \quad q_1\in\mathbb P_{k-2}(K),\quad q_2\in\overline{\mathbb P}_{k-1}(K)\oplus\overline{\mathbb P}_k(K).$$
The norm equivalence on polynomial space implies
\begin{equation}\label{eq:JLowerBound2}
\|\Delta v_h\|_K\sim\|q_1\|_K+\|q_2\|_K.
\end{equation}
For the proof, we first use the standard scaling argument to obtain the norm equivalence on each simplex $T\in\mathcal T_K$ and then sum up them, where the hidden constants are independent of the element size.

Recalling the definition \eqref{2DCVEspace} of space $V_k^{\mathrm c}(K)$, equation \eqref{eq:PiIdentity} and the definition of $\widetilde J_k^K$, we obtain
\begin{align*}
(v_h-P_K(\widetilde J_k^K v_h),\Delta v_h)_K&=(\widetilde J_k^Kv_h-P_K(\widetilde J_k^K v_h),q_1)_K+(\Pi_k^K(v_h-P_K(\widetilde J_k^Kv_h)),q_2)_K\\
&=(\widetilde J_k^Kv_h-P_K(\widetilde J_k^K v_h),q_1)_K+(\Pi_k^K(\widetilde J_k^Kv_h-P_K(\widetilde J_k^Kv_h)),q_2)_K,
\end{align*}
which, together with the boundedness \eqref{PiBound} of $\Pi_k^K$ and the norm equivalence \eqref{eq:JLowerBound2}, leads to
\begin{equation}\label{eq:JLowerBound3}
(v_h-P_K(\widetilde J_k^K v_h),\Delta v_h)_K\lesssim \|\widetilde J_k^Kv_h-P_K(\widetilde J_k^K v_h)\|_K\|\Delta v_h\|_K\lesssim h_K|\widetilde J_k^Kv_h|_{1,K}\|\Delta v_h\|_K.
\end{equation}
Substituting the above inequality into \eqref{eq:JLowerBound1}, we obtain
$$|v_h|_{1,K}^2 \lesssim |\widetilde J_k^Kv_h|_{1,K}| v_h|_{1,K}+h_K|\widetilde J_k^Kv_h|_{1,K}\|\Delta v_h\|_K\lesssim |\widetilde J_k^Kv_h|_{1,K}| v_h|_{1,K},$$
where we used the inverse inequality \eqref{VEInverseInequality}. So the lower bound \eqref{eq:WCFEILowerBound} holds for $m=1$.

For $m=0$, we use the lower bound \eqref{eq:WCFEILowerBound} with $m=1$ and the inverse inequality on polynomials to obtain
\begin{align}\label{eq:JLowerBound4}
\|v_h\|_K &\leq \|v_h-P_Kv_h\|_K+\|P_K(\widetilde J_k^Kv_h)\|_K \nonumber\\
&\lesssim h_K|v_h|_{1,K}+\|\widetilde J_k^Kv_h\|_K\lesssim h_K|\widetilde J_k^Kv_h|_{1,K}+\|\widetilde J_k^Kv_h\|_K\lesssim\|\widetilde J_k^Kv_h\|_K .
\end{align}
The proof is complete.
\end{proof}

\subsection{Interpolation from conforming VE to FE in 3D}

Let $K$ is a polyhedron with the subdivision $\mathcal T_K$. Given a tetrahedron $T\in \mathcal T_K$, we define the conforming FE by setting the local space
$$\widetilde W_k(T)=\mathbb P_k(T)\oplus b_T\overline{\mathbb P}_{k-3}(T)\oplus b_T\overline{\mathbb P}_{k-2}(T)\oplus\sum_{F\in\mathcal F(T)}b_F\overline{\mathbb P}_{k-2}(F)$$
with the DoFs
\begin{align}
&v(\delta),\quad\delta\in\mathcal V(T), \label{W3DCFEDoF1}\\
& \frac{1}{|e|}(v,q)_e,\quad q\in\mathbb M_{k-2}(e),\ e\in\mathcal E(T), \label{W3DCFEDoF2}\\
& \frac{1}{|F|}(v,q)_F,\quad q\in\mathbb M_{k-2}(F),\ F\in\mathcal F(T), \label{W3DCFEDoF3}\\
&\frac{1}{|T|}(v,q)_T,\quad q\in\mathbb M_{k-2}(T).\label{W3DCVEDoF4}
\end{align}
On the subdivision $\mathcal T_K$ of $K$, we introduce the macro Lagrange FE by setting the space
$$\widetilde W_k(\mathcal T_K)=\{v_h\in C^0(K);\ v_h|_T\in \widetilde W_k(T),\ T\in\mathcal T_K\}$$
with the above DoFs.
Then we define the interpolation $\widetilde J_k^K$ from $V_k^{\mathrm c}(K)$ to $\widetilde W_k(\mathcal T_K)$ by letting $\widetilde J_k^Kv_h\in\widetilde W_k(\mathcal T_K)$ satisfy
\begin{align}
&\widetilde J_k^K v_h(\delta)=v_h(\delta),\quad\delta\in\mathcal V (K), \label{Jh3DCLFEIntDoF1}\\
&\widetilde J_k^K v_h(\delta)=\Pi_{k,h}^{\partial K}v_h(\delta),\quad\delta\in\mathcal V^{\partial} (\mathcal T_K)\backslash\mathcal V (K),\\
&\widetilde J_k^K v_h(\delta)=\Pi_k^Kv_h(\delta),\quad\delta\in\mathcal V^{\mathrm{int}} (\mathcal T_K),\label{Jh3DCFEIntDoF3}\\
& \frac{1}{| e|}(\widetilde J_k^K v_h-v_h,q)_e=0,\quad q\in\mathbb M_{k-2}(e),\  e\in\mathcal E(K),\label{Jh3DCFEIntDoF4}\\
& \frac{1}{| e|}(\widetilde J_k^K v_h-\Pi_{k,h}^{\partial K}v_h,q)_e=0,\quad q\in\mathbb M_{k-2}(e),\  e\in\mathcal E^\partial(\mathcal T_K)\backslash \mathcal E(K),\label{Jh3DCFEIntDoF5}\\
& \frac{1}{| e|}(\widetilde J_k^K v_h-\Pi_k^Kv_h,q)_e=0,\quad q\in\mathbb M_{k-2}(e),\  e\in\mathcal E^{\mathrm{int}}(\mathcal T_K),\label{Jh3DCFEIntDoF6}\\
& \frac{1}{| F|}(\widetilde J_k^K v_h-Q_{k,h}^{\partial K}v_h,q)_{ F}=0,\quad q\in\mathbb M_{k-2}(  F),\  F\in\mathcal F^\partial( \mathcal T_K),\label{Jh3DCFEIntDoF7}\\
& \frac{1}{|F|}(\widetilde J_k^K v_h-\Pi_k^Kv_h,q)_{F}=0,\quad q\in\mathbb M_{k-2}( F),\  F\in\mathcal F^{\mathrm{int}}( \mathcal T_K),\\
&\frac{1}{|T|}(\widetilde J_k^Kv_h-Q_k^K v_h,q)_T=0,\quad q\in\mathbb M_{k-2}(T),\ T\in\mathcal T_K.\label{Jh3DCFEIntDoF9}
\end{align}

Similar to the previous discussion in Section \ref{sec:Interpolation3D}, we easily verify the computability and preservation of $\widetilde J_k^K$, i.e. $\widetilde J_k^K$ satisfies the conditions {\bf B1-B2}. By the interpolation theory of FEM, the same discussion as in the proof of Lemma \ref{lem:3DCFEIUpperBound} yields the upper bound
\begin{equation}\label{eq:3DWCFEIUpperBound}
|\widetilde J_k^Kv_h|_{m,K}\lesssim|v_h|_{m,K},\quad \forall v_h\in V_k^{\mathrm c}(K),\ m=0,1.
\end{equation}
Before proving the lower bound of $\widetilde J_k^Kv_h$, we need to present an upper bound of $\Pi_k^Kv_h$.
\begin{lemma}
For $\Pi_k^K$ it holds
\begin{equation}\label{eq:3DPiEstimate}
|\Pi_k^Kv_h|_{m,K}\lesssim |\widetilde J_k^Kv_h|_{m,K},\quad \forall v_h\in V_k^{\mathrm c}(K),\ m=0,1.
\end{equation}
\end{lemma}

\begin{proof}
We set $\widetilde J_k^{\partial K}|_F=\widetilde J_k^F$ for $F\in\mathcal F(K)$.
Similar to the argument in the proof of Lemma \ref{lem:3DinterLowBound}, we have 
\begin{equation}\label{eq:3DPiBound1}
\widetilde J_k^Kv_h|_{\partial K}=\widetilde J_k^{\partial K}v_h,\quad
|v_h|_{m,\partial K}\lesssim |\widetilde J_k^{\partial K}v_h|_{m,\partial K},\quad m=0,1.
\end{equation}
Let $P_Kv_h$ still denote the average values of function $v_h$ at all vertices of $K$ such that the properties \eqref{eq:JLowerBound0} hold. Then the definition of $\widetilde J_k^K$ leads to
\begin{align*}
|\Pi_k^Kv_h|_{1,K}^2 &=(\nabla(v_h-P_Kv_h),\nabla\Pi_k^Kv_h)_K=-(v_h-P_Kv_h,\Delta\Pi_k^K v_h)_K+\int_{\partial K}(v_h-P_Kv_h)\frac{\partial \Pi_k^Kv_h}{\partial\bm n_K}\mathrm ds \\
&=-(\widetilde J_k^Kv_h-P_K(\widetilde J_k^Kv_h),\Delta\Pi_k^K v_h)_K+\int_{\partial K}(v_h-P_Kv_h)\frac{\partial\Pi_k^K v_h}{\partial\bm n_K}\mathrm ds, 
\end{align*}
which, together with the inverse inequality, trace inequality, \eqref{eq:3DPiBound1} and \eqref{eq:JLowerBound0}, yields
\begin{align*}
|\Pi_k^Kv_h|_{1,K}^2 &\leq \|\widetilde J_k^Kv_h-P_K(\widetilde J_k^Kv_h)\|_K\|\Delta\Pi_k^K v_h\|_K+\|v_h-P_Kv_h\|_{\partial K}\Big\|\frac{\partial\Pi_k^K v_h}{\partial\bm n_K}\Big\|_{\partial K}\\
&\lesssim (|\widetilde J_k^Kv_h|_{1,K}+h_K^{-\frac{1}2}\|v_h-P_Kv_h\|_{\partial K})|\Pi_k^Kv_h|_{1,K}\\
&\lesssim |\widetilde J_k^Kv_h|_{1,K}|\Pi_k^Kv_h|_{1,K}.
\end{align*}
This leads to the estimate \eqref{eq:3DPiEstimate} for $m=1$.

For $m=0$, recalling the fact that $Q_0^{\partial K}v_h=Q_0^{\partial K}(\Pi_k^Kv_h)$, we use the Poincar\'e-Friedrichs inequality, inverse inequality and \eqref{eq:3DPiBound1} to obtain
\begin{align*}
\|\Pi_k^Kv_h\|_K &\leq \|\Pi_k^Kv_h-Q_0^{\partial K}(\Pi_k^Kv_h)\|_K+\|Q_0^{\partial K}v_h\|_K\\
&\lesssim h_K|\Pi_k^Kv_h|_{1,K}+h_K^{\frac{1}2}\|v_h\|_{\partial K} \\
&\lesssim h_K|\Pi_k^Kv_h|_{1,K}+\|\widetilde J_k^Kv_h\|_K,
\end{align*}
which, together with the estimate \eqref{eq:3DPiEstimate} for $m=1$ and the inverse inequality, concludes the proof.
\end{proof}

Next we derive the lower bound of $\widetilde J_k^K$.

\begin{lemma}
For any given $v_h \in V_k^{\mathrm c}(K)$, it holds
\begin{equation}\label{eq:3DWCFEILowerBound}
|v_h|_{m,K}\lesssim|\widetilde J_k^Kv_h|_{m,K},\quad m=0,1.
\end{equation}
\end{lemma}

\begin{proof}
First, we have
\begin{equation}\label{eq:3DJhLowBound1}
|v_h|_{1,K}^2=-(v_h-P_K(\widetilde J_k^Kv_h),\Delta v_h)_K+\int_{\partial K}(v_h-P_Kv_h)\frac{\partial v_h}{\partial\bm n_K}\mathrm ds.
\end{equation}
For the first term in \eqref{eq:3DJhLowBound1}, we use the similar argument in proof of estimate \eqref{eq:JLowerBound3} to obtain
$$(v_h-P_K(\widetilde J_k^Kv_h),\Delta v_h)_K\lesssim (|\widetilde J_k^Kv_h|_{1,K}+h_K^{-1}\|\Pi_k^K(v_h-P_K v_h)\|_K)|v_h|_{1,K}$$
which, combining the estimates \eqref{eq:3DPiEstimate} with $m=1$, yields
\begin{equation}\label{eq:3DJhLowBound2}
(v_h-P_K(\widetilde J_k^Kv_h),\Delta v_h)_K\lesssim|\widetilde J_k^Kv_h|_{1,K}|v_h|_{1,K}.
\end{equation}
For the second term in \eqref{eq:3DJhLowBound1},  similar to the proof of estimate \eqref{eq:3DGreenFormula3}, we obtain
\begin{equation}\label{eq:3DJhLowBound3}
\int_{\partial K}(v_h-P_Kv_h)\frac{\partial v_h}{\partial\bm n_K}\mathrm ds\lesssim |\widetilde J_k^Kv_h|_{1,K}|v_h|_{1,K}.
\end{equation}
Substituting the estimates \eqref{eq:3DJhLowBound2}-\eqref{eq:3DJhLowBound3} into \eqref{eq:3DJhLowBound1}, we get the lower bound \eqref{eq:3DWCFEILowerBound} for $m=1$.

For $m=0$, the same discussion as in the proof of \eqref{eq:JLowerBound4} yields the lower bound \eqref{eq:3DWCFEILowerBound}.
\end{proof}

From the inequalities \eqref{eq:3DWCFEIUpperBound} and \eqref{eq:3DWCFEILowerBound}, we see that the interpolation $\widetilde J_k^K$ satisfies the condition {\bf B3}.

\subsection{Interpolation from nonconforming VE to FE}

Let $K$ be a polygon for $d=2$ and a polyhedron for $d=3$. We use the macro FE space $\widetilde W_{k+1}(\mathcal T_K)$ of degree $(k+1)$ as the interpolation space to define the interpolation from the nonconforming VE to FE. For a given VE function $v_h\in V_k^{\mathrm{nc}}(K)$, we define the interpolation $\widetilde J_k^{\mathrm{nc},K}$ from $V_k^{\mathrm{nc}}(K)$ to $\widetilde W_{k+1}(\mathcal T_K)$ by letting $\widetilde J_k^{\mathrm{nc},K}v_h\in\widetilde W_{k+1}(\mathcal T_K)$ satisfy: for $d=2$
\begin{align*}
&\widetilde J_k^{\mathrm{nc},K} v_h(\delta)=v_h(\delta),\quad\delta\in\mathcal V(K), \\
&\widetilde J_k^{\mathrm{nc},K} v_h(\delta)=\Pi_k^Kv_h(\delta),\quad\delta\in\mathcal V^{\mathrm{int}}(\mathcal T_K),\\
& \frac{1}{|e|}(\widetilde J_k^{\mathrm{nc},K} v_h-v_h,q)_e=0,\quad q\in\mathbb M_{k-1}(e),\ e\in\mathcal E^\partial(\mathcal T_K), \\
& \frac{1}{|e|}(\widetilde J_k^{\mathrm{nc},K} v_h-\Pi_k^Kv_h,q)_e=0,\quad q\in\mathbb M_{k-1}(e),\ e\in\mathcal E^{\mathrm{int}}(\mathcal T_K), \\
&\frac{1}{|T|}(\widetilde J_k^{\mathrm{nc},K}v_h-Q_{k+1}^K v_h,q)_T=0,\quad q\in\mathbb M_{k-1}(T),\ T\in\mathcal T_K,
\end{align*}
and for $d=3$
\begin{align*}
&\widetilde J_k^{\mathrm{nc},K} v_h(\delta)=v_h(\delta),\quad\delta\in\mathcal V(K), \\
&\widetilde J_k^{\mathrm{nc},K} v_h(\delta)=\Pi_{k+1,h}^{\partial K} v_h(\delta),\quad\delta\in\mathcal V^\partial(\mathcal T_K)\backslash\mathcal V(K), \\
&\widetilde J_k^{\mathrm{nc},K} v_h(\delta)=\Pi_k^K v_h(\delta),\quad\delta\in\mathcal V^{\mathrm{int}}(\mathcal T_K), \\
& \frac{1}{|e|}(\widetilde J_k^{\mathrm{nc},K} v_h-v_h,q)_e=0,\quad q\in\mathbb M_{k-1}(e),\ e\in\mathcal E(K), \\
& \frac{1}{|e|}(\widetilde J_k^{\mathrm{nc},K} v_h-\Pi_{k+1,h}^{\partial K}v_h,q)_e=0,\quad q\in\mathbb M_{k-1}(e),\ e\in\mathcal E^\partial(\mathcal T_K)\backslash\mathcal E( K), \\
& \frac{1}{|e|}(\widetilde J_k^{\mathrm{nc},K} v_h-\Pi_k^Kv_h,q)_e=0,\quad q\in\mathbb M_{k-1}(e),\ e\in\mathcal E^{\mathrm{int}}(\mathcal T_K),\\
&\frac{1}{| F|}(\widetilde J_k^{\mathrm{nc},K} v_h-Q_{k+1,h}^{\partial K}v_h,q)_{ F}=0,\quad q\in\mathbb M_{k-1}( F),\  F\in\mathcal F^{\partial}(\mathcal T_K),\\
&\frac{1}{|F|}(\widetilde J_k^{\mathrm{nc},K} v_h-\Pi_k^Kv_h,q)_F=0,\quad q\in\mathbb M_{k-1}(F),\ F\in\mathcal F^{\mathrm{int}}(\mathcal T_K),\\
&\frac{1}{|T|}(\widetilde J_k^{\mathrm{nc},K} v_h-Q_{k+1}^K v_h,q)_T=0,\quad q\in\mathbb M_{k-1}(T),\ T\in\mathcal T_K.
\end{align*}

Following the previous discussions, we can prove that the interpolation $\widetilde J_k^{\mathrm{nc},K}$ satisfies the conditions {\bf B1-B3}.

\begin{remark}\label{rmk:PbFEM}
If we set $W_h=\widetilde J_h V_h$, the discrete space $W_h$ has the same unisolvent DoFs as the VE space $V_h$, meanwhile the shape functions in $W_h$ have the explicit expression which are piecewise polynomials on the subdivision of each polytope. The interpolation on $W_h$ can be defined by $J_h=\widetilde J_hI_h$ with the estimates
$$|v-J_hv|_{m,K}\lesssim |v-I_hv|_{m,K}+\inf_{q\in\mathbb P_k(K)}|v-q|_{m,K},\quad K\in\mathcal T_h,\ m=0,1.$$
The second stabilization-free scheme \eqref{eq:SFVEM2} based on $\widetilde J_h$ inspires us to develop a new scheme directly on space $W_h$, which reads: find $u_h\in W_h$ such that
\begin{equation*}
(\alpha\nabla_h u_h, \nabla_h v_h)+(\beta u_h,v_h)=( f,v_h),\quad\forall v_h\in W_h
\end{equation*}
with the following error estimates
\begin{equation*}
\|u-u_h\|_{1,h}\lesssim \|u-J_hu\|_{1,h}+\sup_{v_h\in W_h\backslash\{0\}}\left(\frac{1}{\|v_h\|_{1,h}}\sum_{K\in\mathcal T_h}\int_{\partial K}\alpha\frac{\partial u}{\partial\bm n_K} v_h\mathrm ds\right),
\end{equation*}
which is the same as that in the classic FEM on simplicial meshes (see C\'ea's lemma for the conforming case and second Strang lemma for the nonconforming case in \cite{Ciarlet1978}).
We observe that the construction of $\widetilde J_h$ mainly depends on the projections $\Pi_k$ and $Q_k$, so we call this new method as {\bf the projection-based FEM} on polytopal meshes.

\end{remark}

\section{Numerical tests}
\label{sec:NumerTest}

In this section, we present numerical experiments to confirm the theoretical convergence rates of the proposed stabilization-free conforming and nonconforming VEMs. We also investigate the performance of the proposed method through comparisons with the standard VEM for an anisotropic elliptic problem and the standard FEM for an internal-layer problem.

\subsection{Convergence tests}
The second-order elliptic problem \eqref{eq:2ndEquation} with $\alpha=1$ and $\beta=1$ is solved by the stabilization-free conforming and nonconforming VEMs with $k=1,2,3$ on the unit domain $\Omega=(0,1)^d$. The right-hand side $f$ is chosen so that the exact solution is
\begin{equation*} 
	u(\bm x) =
	\begin{cases}
		\sin(\pi x_1)\sin(\pi x_2), & d=2,\\
		\sin(\pi x_1)\sin(\pi x_2)\sin(\pi x_3), & d=3.
	\end{cases}
\end{equation*}
Note that for the first type of stabilization-free VEM \eqref{eq:SFVEM1} with consistency, the numerical results are similar. Therefore, for brevity, here we only present the numerical results for the second type of the stabilization-free VEM \eqref{eq:SFVEM2} without consistency.

\begin{figure}[ht]
	\centering
	\includegraphics[width=0.45\textwidth]{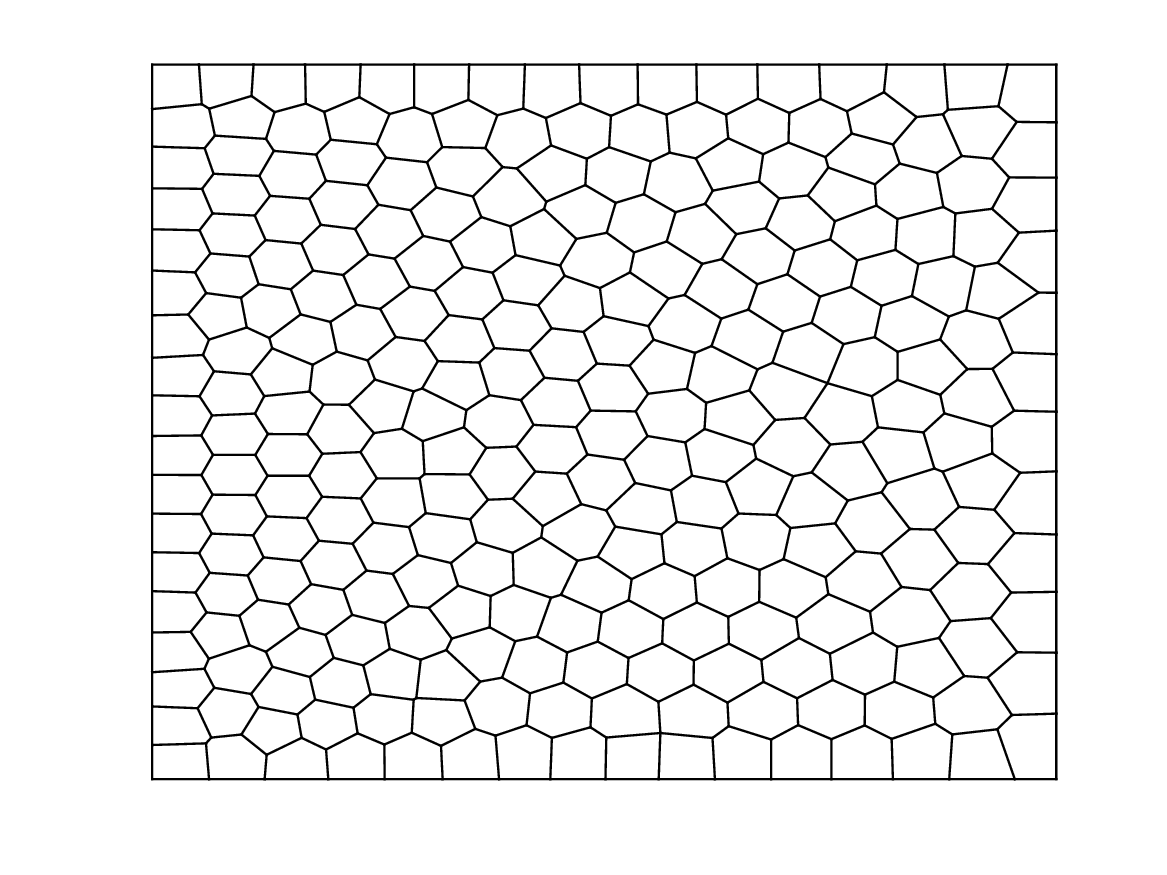}
	\includegraphics[width=0.45\textwidth]{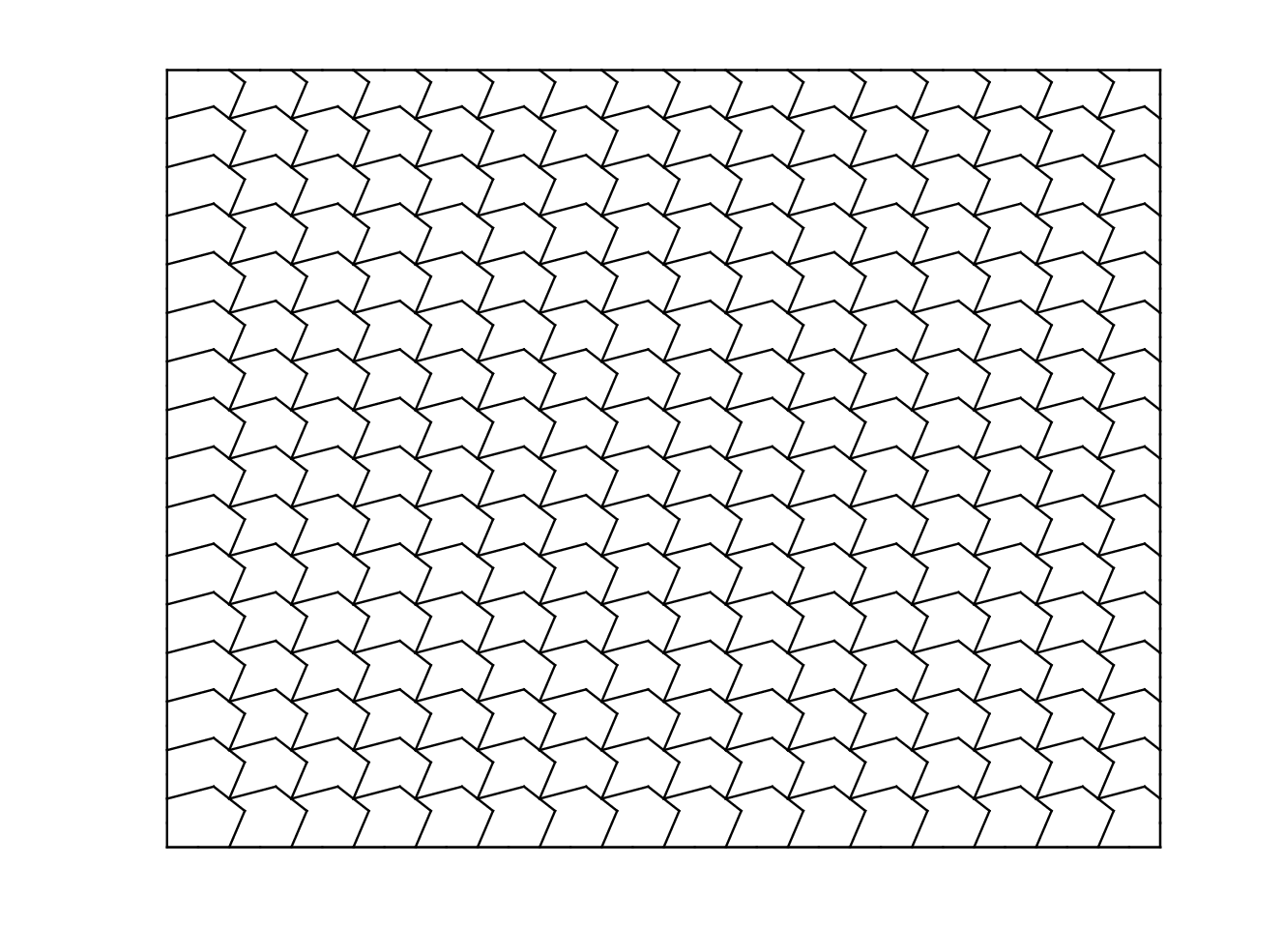}
	\caption{The Voronoi (left) and nonconvex (right) meshes.}
	\label{figure=UnipolyMesh=NonConvexMesh}
\end{figure}

For the two-dimensional case, the unit square $\Omega$ is partitioned into the Voronoi mesh $\mathcal{T}_{h,\mathrm{2D}}^1$ and the nonconvex polygonal mesh $\mathcal{T}_{h,\mathrm{2D}}^2$, as shown in Figure \ref{figure=UnipolyMesh=NonConvexMesh}. The simplicial submesh of each polygon is constructed by connecting the barycenter of each polygon to its vertices.\\

\begin{figure}[ht]
	\centering
	\includegraphics[width=0.35\textwidth]{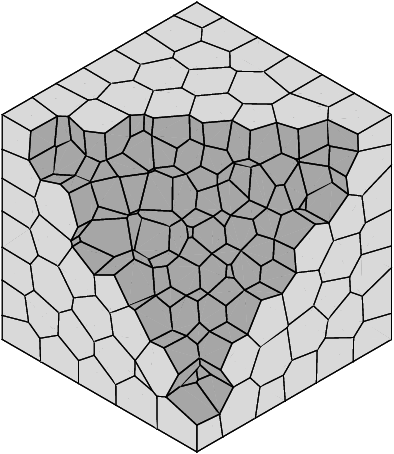}\hspace{3em}
	\includegraphics[width=0.35\textwidth]{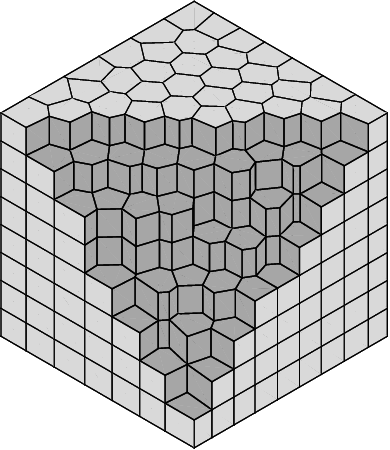}
	\caption{The CVT (left) and prism (right) meshes.}
	\label{figure3D=CVTMesh=PrismMesh}
\end{figure}

\begin{figure}[ht]
	\centering
	\begin{subfigure}[t]{0.49\textwidth} 
		\centering
		\includegraphics[width=0.9\textwidth]{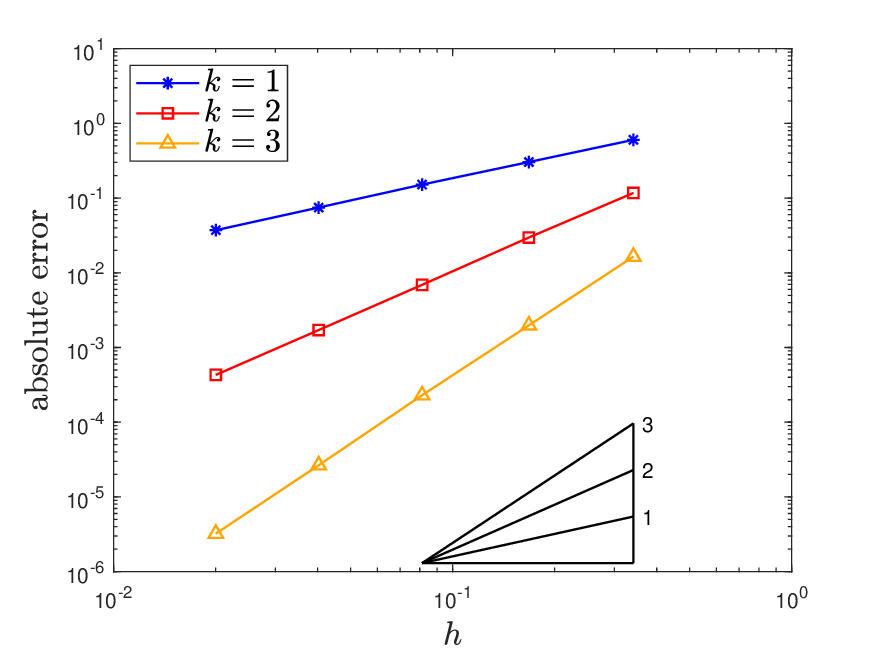}
		\caption{Voronoi Mesh $\mathcal{T}_{h,\mathrm{2D}}^1$}
	\end{subfigure}
	\hfill
	\begin{subfigure}[t]{0.49\textwidth}
		\centering
		\includegraphics[width=0.9\textwidth]{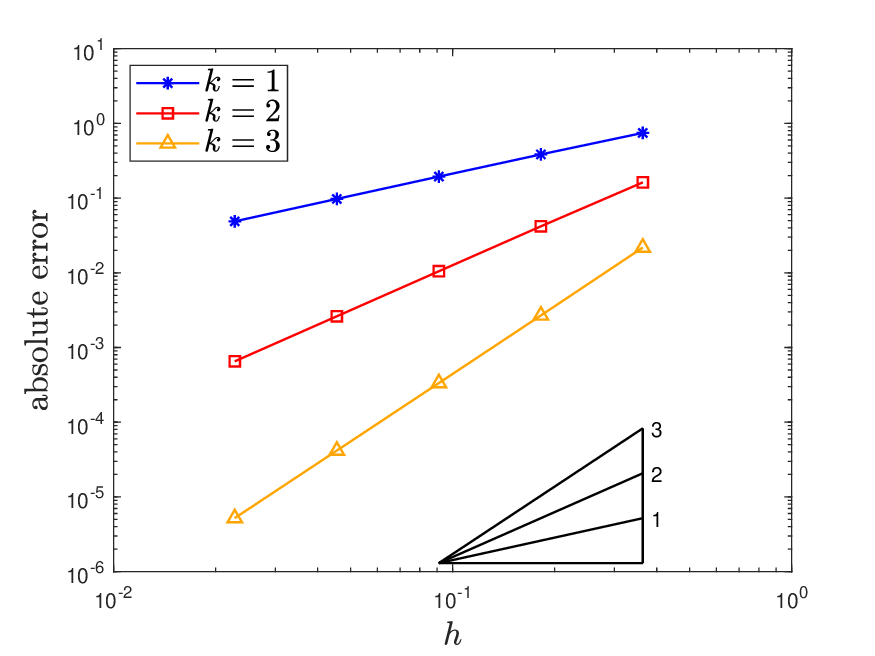}
		\caption{Nonconvex Mesh $\mathcal{T}_{h,\mathrm{2D}}^2$}
	\end{subfigure}

	\begin{subfigure}[t]{0.49\textwidth}
		\centering
		\includegraphics[width=0.9\textwidth]{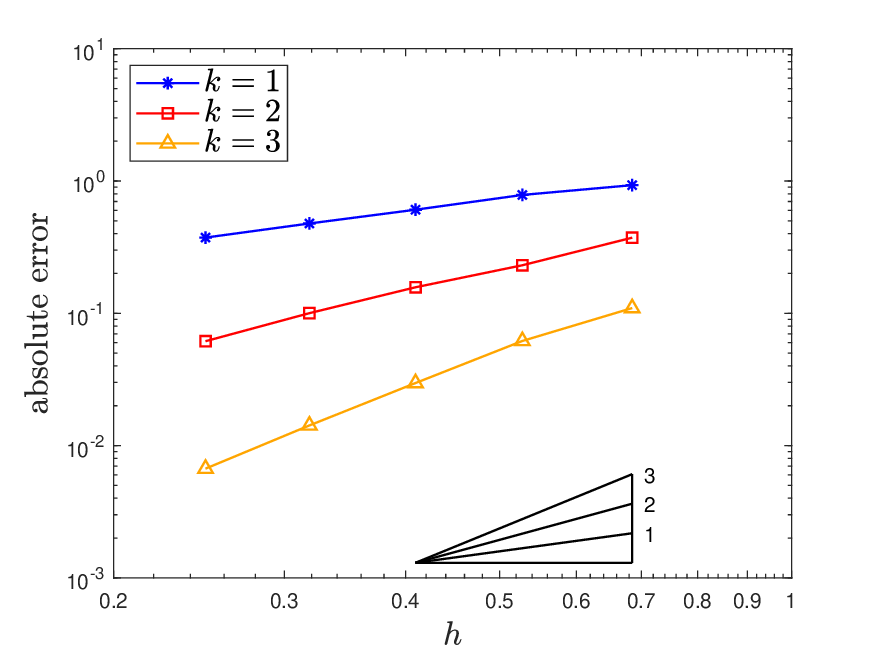}
		\caption{CVT Mesh $\mathcal{T}_{h,\mathrm{3D}}^1$}
	\end{subfigure}
	\hfill
	\begin{subfigure}[t]{0.49\textwidth}
		\centering
		\includegraphics[width=0.9\textwidth]{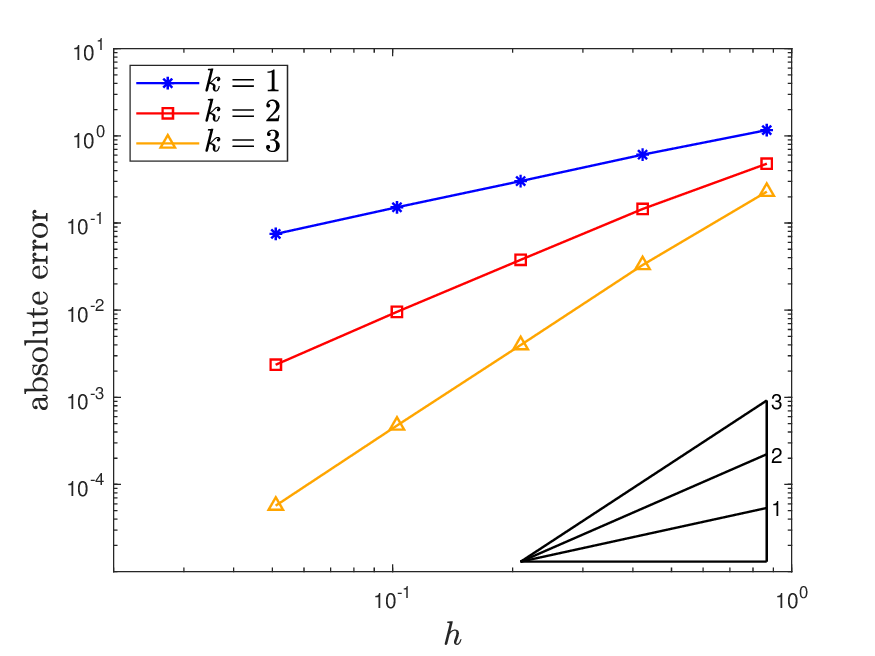}
		\caption{Prism Mesh $\mathcal{T}_{h,\mathrm{3D}}^2$}
	\end{subfigure}
	\caption{The $H^1$ error $|{u}-\widetilde J_h{u}_h|_{1}$ of the stabilization-free conforming VEM without consistency with $k=1,2,3$.}
	\label{Figure=H1Error}
\end{figure}

\begin{figure}[ht]
	\centering
	\begin{subfigure}[t]{0.49\textwidth}
		\centering
		\includegraphics[width=0.9\textwidth]{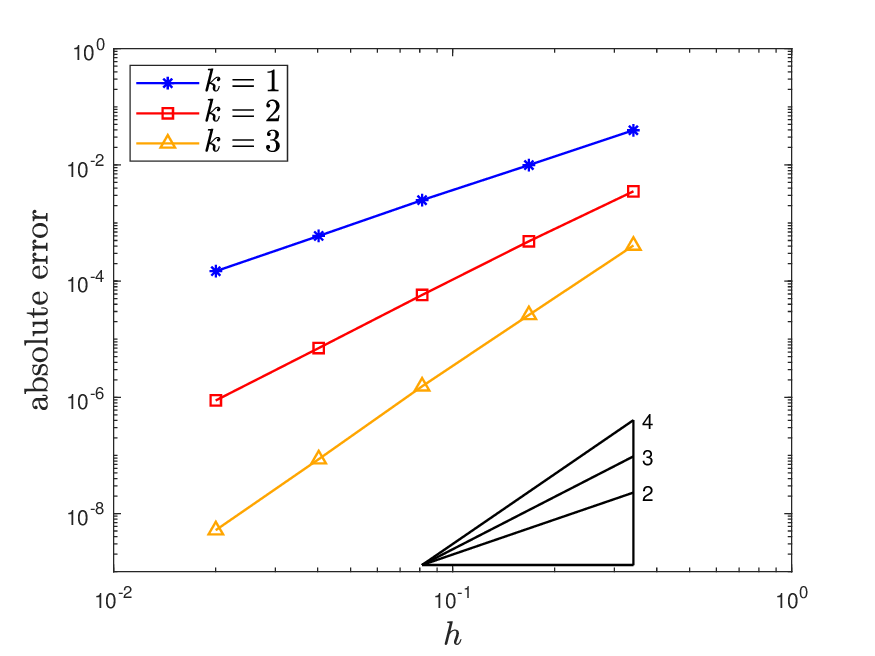}
		\caption{Voronoi Mesh $\mathcal{T}_{h,\mathrm{2D}}^1$}
	\end{subfigure}
	\hfill
	\begin{subfigure}[t]{0.49\textwidth}
		\centering
		\includegraphics[width=0.9\textwidth]{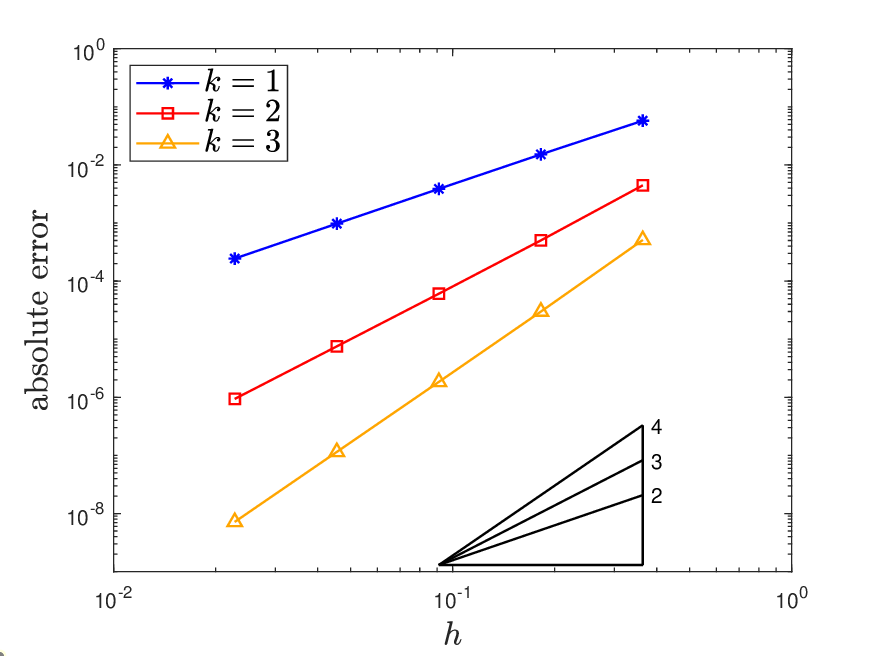}
		\caption{Nonconvex Mesh $\mathcal{T}_{h,\mathrm{2D}}^2$}
	\end{subfigure}

	\begin{subfigure}[t]{0.49\textwidth}
		\centering
		\includegraphics[width=0.9\textwidth]{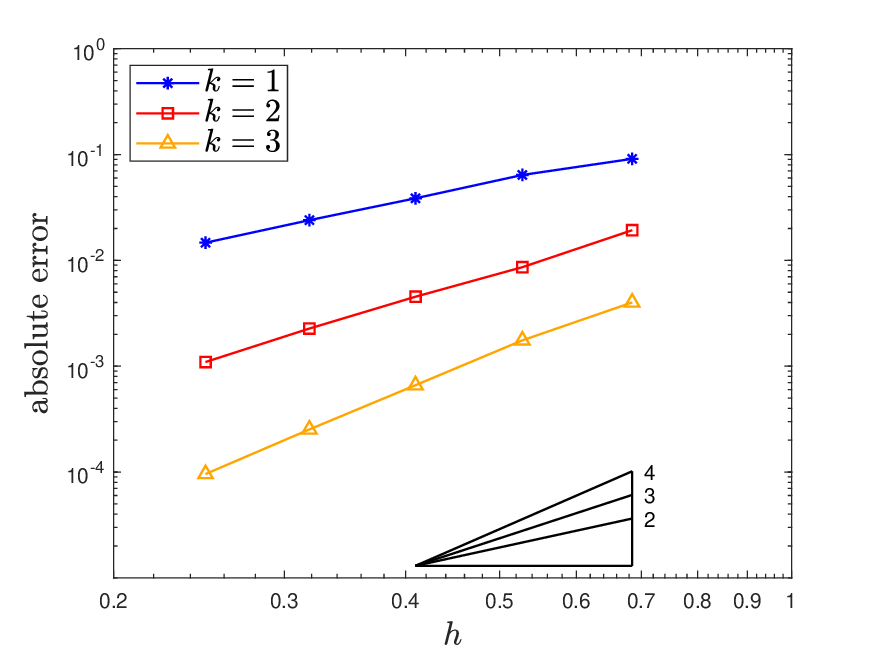}
		\caption{CVT Mesh $\mathcal{T}_{h,\mathrm{3D}}^1$}
	\end{subfigure}
	\hfill
	\begin{subfigure}[t]{0.49\textwidth}
		\centering
		\includegraphics[width=0.9\textwidth]{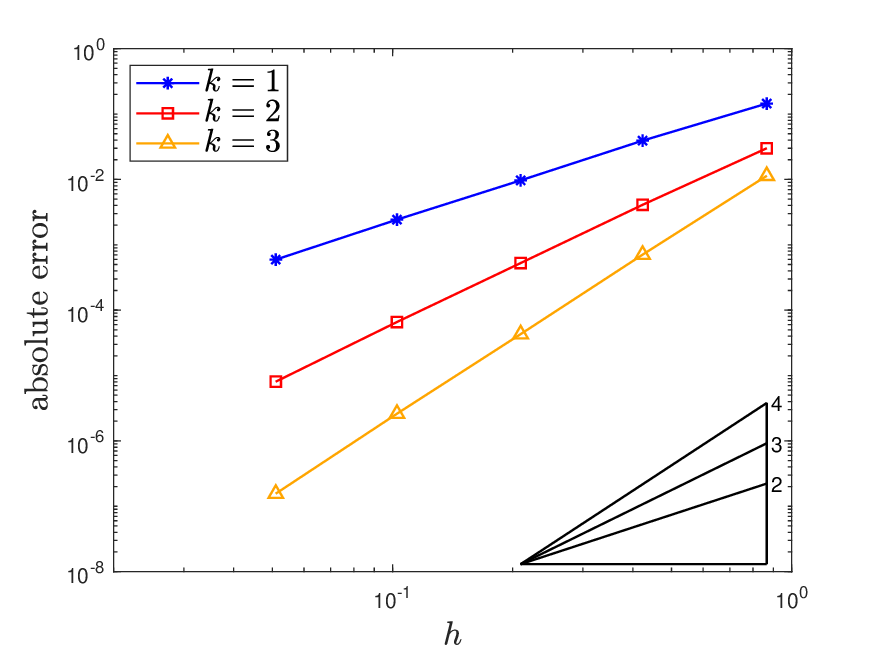}
		\caption{Prism Mesh $\mathcal{T}_{h,\mathrm{3D}}^2$}
	\end{subfigure}
	\caption{The $L^2$ error $\| {u}-\widetilde J_h{u}_h \|$ of the stabilization-free conforming VEM without consistency with $k=1,2,3$.}
	\label{Figure=H0Error}
\end{figure}

\vspace{-0.3cm}

\begin{figure}[ht] 
	\centering
	\begin{subfigure}[t]{0.49\textwidth}
		\centering
		\includegraphics[width=0.9\textwidth]{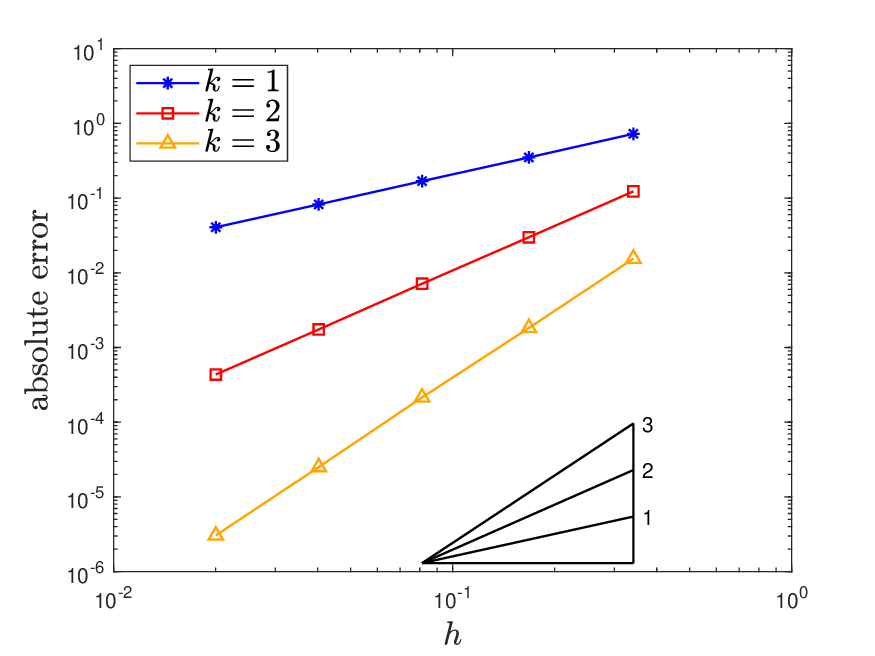}
		\caption{Voronoi Mesh $\mathcal{T}_{h,\mathrm{2D}}^1$}
	\end{subfigure}
	\hfill
	\begin{subfigure}[t]{0.49\textwidth}
		\centering
		\includegraphics[width=0.9\textwidth]{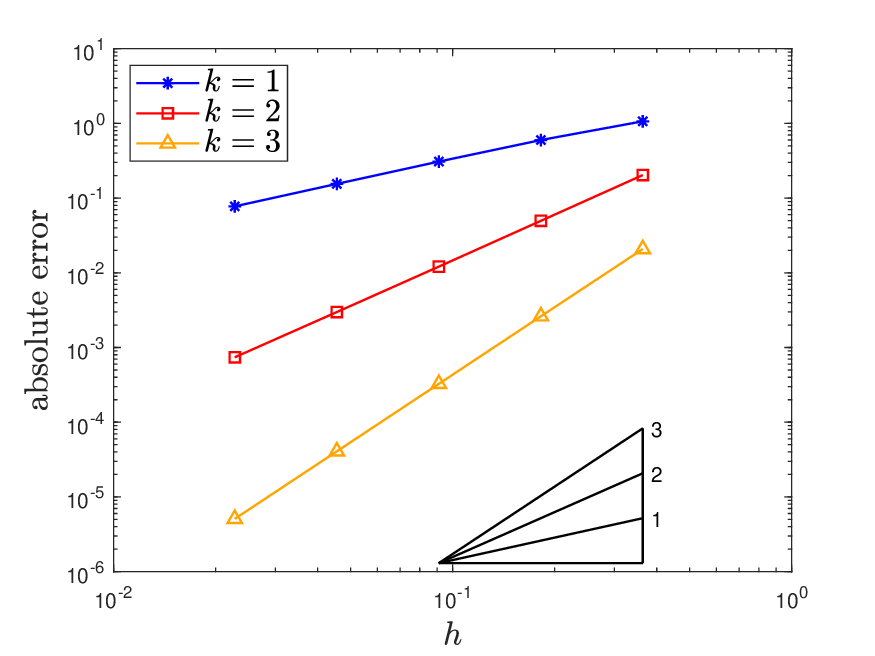}
		\caption{Nonconvex Mesh $\mathcal{T}_{h,\mathrm{2D}}^2$}
	\end{subfigure}

	\begin{subfigure}[t]{0.49\textwidth}
		\centering
		\includegraphics[width=0.9\textwidth]{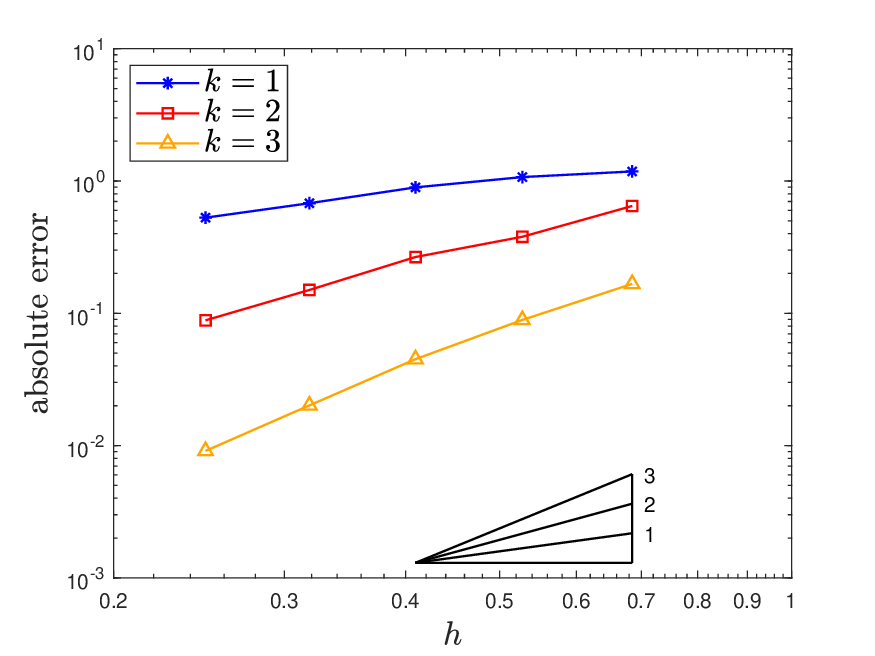}
		\caption{CVT Mesh $\mathcal{T}_{h,\mathrm{3D}}^1$}
	\end{subfigure}
	\hfill
	\begin{subfigure}[t]{0.49\textwidth}
		\centering
		\includegraphics[width=0.9\textwidth]{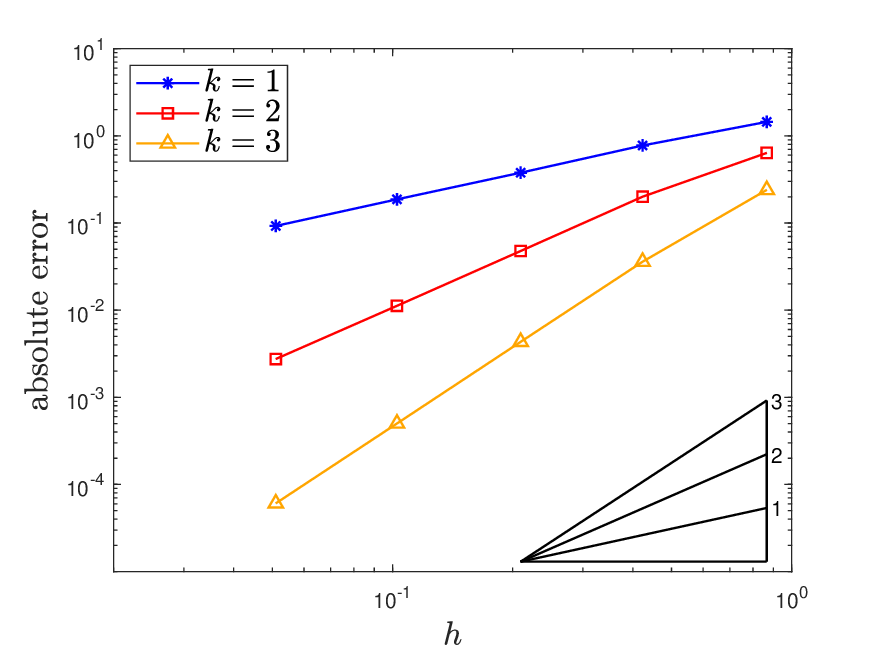}
		\caption{Prism Mesh $\mathcal{T}_{h,\mathrm{3D}}^2$}
	\end{subfigure}
	\caption{The $H^1$ error $|{u}-\widetilde J_h{u}_h|_{1}$ of the stabilization-free nonconforming VEM without consistency with $k=1,2,3$.}
	\label{Figure=NCH1Error}
\end{figure}

For the three-dimensional case, the unit cube $\Omega$ is partitioned into the CVT mesh $\mathcal{T}_{h,\mathrm{3D}}^1$ and the prism mesh $\mathcal{T}_{h,\mathrm{3D}}^2$, as shown in Figure \ref{figure3D=CVTMesh=PrismMesh}. For each polyhedron, the simplicial submesh is constructed by first connecting the barycenter of each polygonal face to its vertices, thereby subdividing the face into triangles, and then connecting the barycenter of the polyhedron to the vertices of each resulting triangle to form tetrahedra.

\begin{figure}[ht]
	\centering
	\begin{subfigure}[t]{0.48\textwidth}
		\centering
		\includegraphics[width=0.87\textwidth]{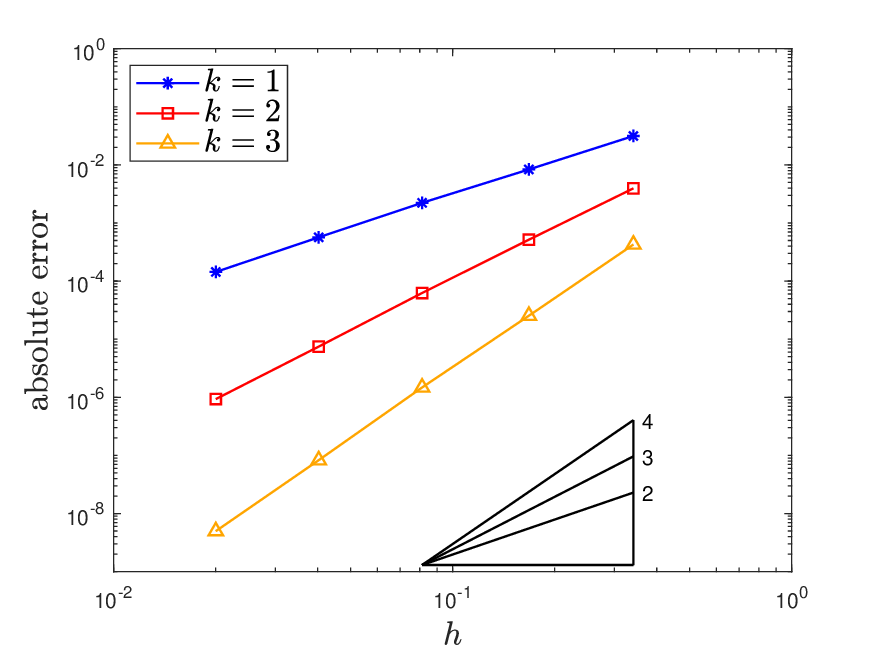}
		\caption{Voronoi Mesh $\mathcal{T}_{h,\mathrm{2D}}^1$}
	\end{subfigure}
	\hfill
	\begin{subfigure}[t]{0.48\textwidth}
		\centering
		\includegraphics[width=0.87\textwidth]{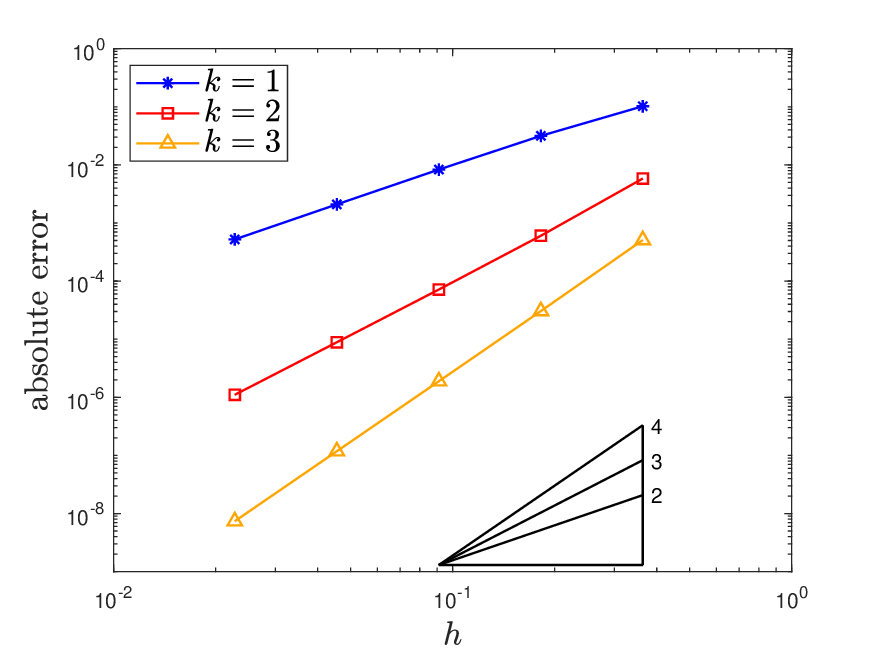}
		\caption{Nonconvex Mesh $\mathcal{T}_{h,\mathrm{2D}}^2$}
	\end{subfigure}
    
    \vspace{0.3cm}

	\begin{subfigure}[t]{0.48\textwidth}
		\centering
		\includegraphics[width=0.87\textwidth]{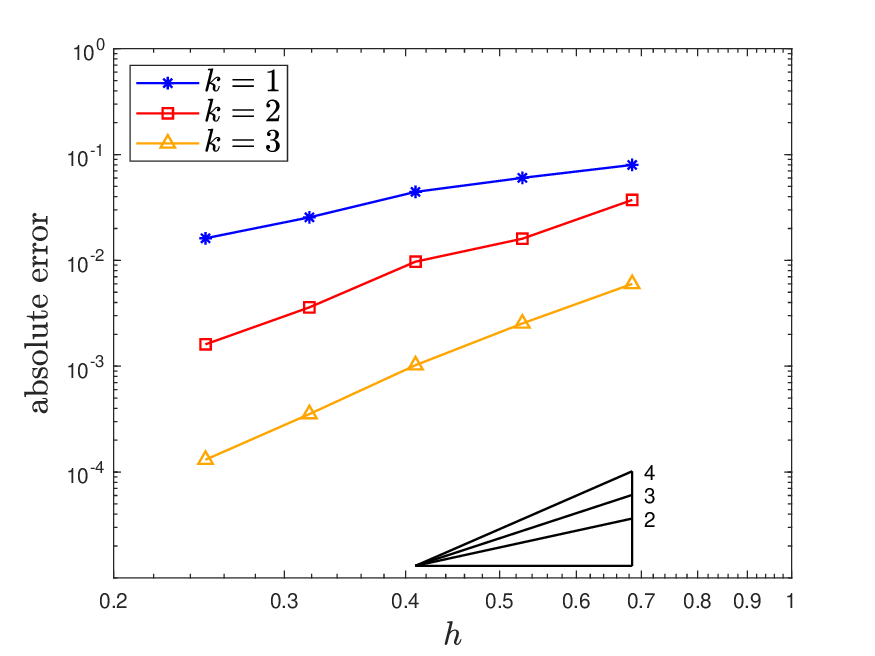}
		\caption{CVT Mesh $\mathcal{T}_{h,\mathrm{3D}}^1$}
	\end{subfigure}
	\hfill
	\begin{subfigure}[t]{0.48\textwidth}
		\centering
		\includegraphics[width=0.87\textwidth]{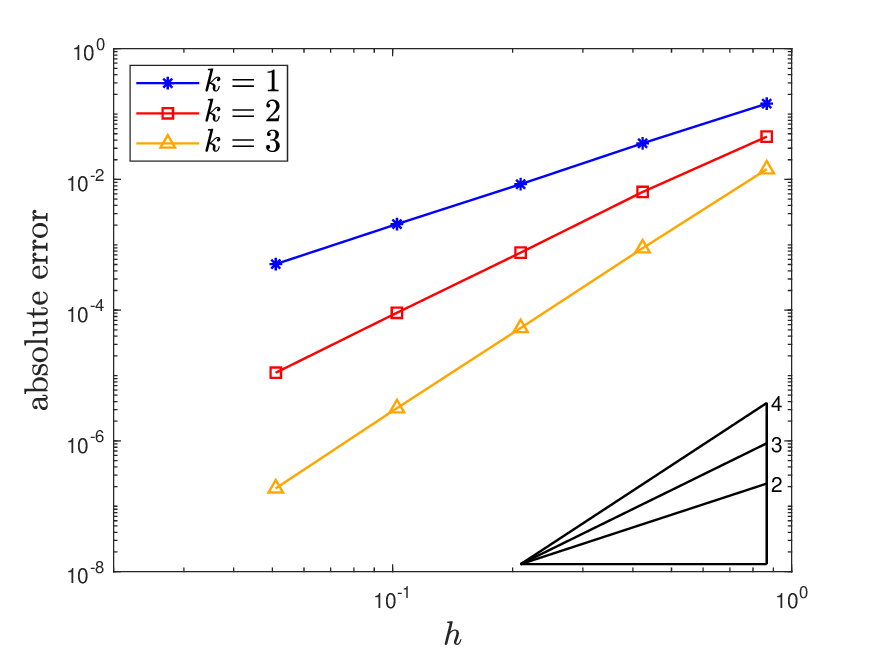}
		\caption{Prism Mesh $\mathcal{T}_{h,\mathrm{3D}}^2$}
	\end{subfigure}
    
	\caption{The $L^2$ error $\| {u}-\widetilde J_h{u}_h \|$ of the stabilization-free nonconforming VEM without consistency with $k=1,2,3$.}
	\label{Figure=NCH0Error}
\end{figure}

Figures~\ref{Figure=H1Error}-\ref{Figure=NCH0Error} present the numerical results for the stabilization-free conforming and nonconforming VEMs \eqref{eq:SFVEM2} without the consistency term on the mesh families $\mathcal{T}_{h,\mathrm{2D}}^1$ and $\mathcal{T}_{h,\mathrm{2D}}^2$ in two dimensions, and $\mathcal{T}_{h,\mathrm{3D}}^1$ and $\mathcal{T}_{h,\mathrm{3D}}^2$ in three dimensions. For ease of verification, we plot the absolute error curves against the mesh size $h$ on a log-log scale, with numbers near the triangles indicating the theoretical convergence rates. From the numerical results, we observe that the $H^1$ error $|{u}-\widetilde J_h{u}_h|_{1}$ and the $L^2$ error $\| {u}-\widetilde J_h{u}_h \|$ are $\mathcal{O}(h^k)$ and $\mathcal{O}(h^{k+1})$, respectively, which are consistent with the theoretical predictions.

\subsection{Numerical comparison of VEM and stabilization-free VEM}
we compare the performance of the conforming stabilization-free VEM (SFVEM) and the standard conforming VEM for the following two-dimensional anisotropic elliptic problem:
\begin{equation*}
	\left\{
	\begin{aligned}
		-\nabla\cdot(\mathcal{K}_{\kappa}\nabla u) &= f
		&& \text{in } \Omega,\\
		u &=0
		&& \text{on } \partial\Omega,
	\end{aligned}
	\right.
\end{equation*}
where $\mathcal{K}_{\kappa}$ is the anisotropic diffusion tensor defined by
\begin{equation*}
	\mathcal{K}_{\kappa}
	=
	\operatorname{diag}(1,\kappa).
\end{equation*}
Here, $\kappa$ is the anisotropy parameter. We consider the cases $\kappa=1,10^2$, and $10^4$, where $\kappa=1$ corresponds to the isotropic case and larger values of $\kappa$ represent stronger anisotropy. The right-hand side $f$ is chosen so that the exact solution is
\begin{equation*}
	u(\bm x)=\sin(\pi x_1)\sin(6\pi x_2).
\end{equation*}

For both methods, we set $k=1$ and use the polygonal mesh family $\mathcal{T}_{h,\mathrm{2D}}^1$  introduced in the previous subsection.
For the standard conforming VEM, the local discrete bilinear forms on a polygon $K$ are given by
\begin{equation*}
	a_K^{\mathrm{VEM}}(v_h,w_h)
	=
	\left(\mathcal{K}_{\kappa} Q_{k-1}^K\nabla v_h,Q_{k-1}^K \nabla w_h\right)_K
	+
	S^K\left(v_h-Q_k^K v_h,w_h-Q_k^K w_h\right),
\end{equation*}
and
\begin{equation*}
	b_K^{\mathrm{VEM}}(v_h,w_h)
	=\left(Q_k^K v_h,Q_k^K w_h\right)_K+h_K^2S^K\left(v_h-Q_k^K v_h,w_h-Q_k^K w_h\right),
\end{equation*}
where $S^K:V_k^c(K)\times V_k^c(K)\to\mathbb{R}$ denotes the local DoF--DoF stabilization term defined by
\begin{equation*}
	S^K(v_h,w_h)=\bm\chi^K(v_h)\cdot\bm\chi^K(w_h),
\end{equation*}
with $\bm\chi^K(v_h)$ denoting the vector of degrees of freedom of VE function $v_h$ on $K$. 
Let $\Pi_h$ and $Q_h$ denote the broken $H^1$- and $L^2$-projection operators onto the piecewise polynomials space of degree at most $k$, respectively. 
The right-hand side is discretized by $(f,Q_h v_h)$. Note that here the stabilization parameter in the stabilization term is taken to be 1 as usual.

We recall that a virtual element function is not explicitly known in the interior of each element, the error of VEM cannot in general be evaluated pointwise. To provide a consistent and fully computable criterion for comparing the two methods, we compute the same projection-based errors 
\begin{equation*} 
E_{H^1}^{\mathrm{M}} = \left|u-\Pi_{h}u_h^{\mathrm{M}}\right|_{1}, \qquad 
E_{L^2}^{\mathrm{M}} = \left\|u-Q_{h}u_h^{\mathrm{M}}\right\|, \qquad \mathrm{M}\in\{\mathrm{VEM},\mathrm{SFVEM}\}
\end{equation*} 
for the standard VEM and the stabilization-free VEM (denoted by SFVEM for short), respectively.
The numerical results are shown in Tables \ref{tab:kappa-1}-\ref{tab:kappa-1e4}, where the convergence rates are also reported.

\begin{table}[ht]
\centering
\caption{Numerical results for the anisotropic elliptic problem with $\mathcal{K}_{\kappa}=\operatorname{diag}(1,1)$.}\label{tab:kappa-1}
\resizebox{\textwidth}{!}{%
\begin{tabular}{c|cccc|cccc}
\hline
$h$ & $E_{H^1}^{\mathrm{VEM}}$ & Rate & $E_{H^1}^{\mathrm{SFVEM}}$ & Rate & $E_{L^2}^{\mathrm{VEM}}$ & Rate & $E_{L^2}^{\mathrm{SFVEM}}$ & Rate\\
\hline
0.34 & $9.32111\mathrm{e}{+0}$ & --   & $9.29206\mathrm{e}{+0}$ & --   & $5.10527\mathrm{e}{-1}$ & --   & $5.05340\mathrm{e}{-1}$ & --\\
0.17 & $6.03338\mathrm{e}{+0}$ & 0.61 & $6.05314\mathrm{e}{+0}$ & 0.60 & $1.78146\mathrm{e}{-1}$ & 1.48 & $1.90852\mathrm{e}{-1}$ & 1.37\\
0.08 & $3.22134\mathrm{e}{+0}$ & 0.87 & $3.22128\mathrm{e}{+0}$ & 0.87 & $4.66203\mathrm{e}{-2}$ & 1.85 & $5.15148\mathrm{e}{-2}$ & 1.81\\
0.04 & $1.60866\mathrm{e}{+0}$ & 0.99 & $1.60825\mathrm{e}{+0}$ & 0.99 & $1.09059\mathrm{e}{-2}$ & 2.07 & $1.24359\mathrm{e}{-2}$ & 2.02\\
0.02 & $8.13782\mathrm{e}{-1}$ & 0.98 & $8.13585\mathrm{e}{-1}$ & 0.98 & $2.80015\mathrm{e}{-3}$ & 1.95 & $3.19224\mathrm{e}{-3}$ & 1.95\\
\hline
\end{tabular}}
\end{table}

\begin{table}[ht]
\centering
\caption{Numerical results for the anisotropic elliptic problem with $\mathcal{K}_{\kappa}=\operatorname{diag}(1,10^{2})$.}\label{tab:kappa-1e2}
\resizebox{\textwidth}{!}{%
\begin{tabular}{c|cccc|cccc}
\hline
$h$ & $E_{H^1}^{\mathrm{VEM}}$ & Rate & $E_{H^1}^{\mathrm{SFVEM}}$ & Rate & $E_{L^2}^{\mathrm{VEM}}$ & Rate & $E_{L^2}^{\mathrm{SFVEM}}$ & Rate\\
\hline
0.34 & $1.21227\mathrm{e}{+1}$ & --      & $9.29887\mathrm{e}{+0}$ & --   & $8.62339\mathrm{e}{-1}$ & --      & $5.30676\mathrm{e}{-1}$ & --\\
0.17 & $3.68443\mathrm{e}{+1}$ & $-1.57$ & $6.08720\mathrm{e}{+0}$ & 0.60 & $2.08272\mathrm{e}{+0}$ & $-1.24$ & $1.64552\mathrm{e}{-1}$ & 1.65\\
0.08 & $1.16407\mathrm{e}{+1}$ & 1.59    & $3.22894\mathrm{e}{+0}$ & 0.87 & $3.67504\mathrm{e}{-1}$ & 2.39    & $4.18218\mathrm{e}{-2}$ & 1.89\\
0.04 & $3.95748\mathrm{e}{+0}$ & 1.54    & $1.61371\mathrm{e}{+0}$ & 0.99 & $5.55702\mathrm{e}{-2}$ & 2.69    & $1.02652\mathrm{e}{-2}$ & 2.00\\
0.02 & $1.70819\mathrm{e}{+0}$ & 1.20    & $8.16732\mathrm{e}{-1}$ & 0.97 & $1.17751\mathrm{e}{-2}$ & 2.22    & $2.74095\mathrm{e}{-3}$ & 1.89\\
\hline
\end{tabular}}
\end{table}

\begin{table}[ht]
\centering
\caption{Numerical results for the anisotropic elliptic problem with $\mathcal{K}_{\kappa}=\operatorname{diag}(1,10^{4})$.}\label{tab:kappa-1e4}
\resizebox{\textwidth}{!}{%
\begin{tabular}{c|cccc|cccc}
\hline
$h$ & $E_{H^1}^{\mathrm{VEM}}$ & Rate & $E_{H^1}^{\mathrm{SFVEM}}$ & Rate & $E_{L^2}^{\mathrm{VEM}}$ & Rate & $E_{L^2}^{\mathrm{SFVEM}}$ & Rate\\
\hline
0.34 & $6.77325\mathrm{e}{+2}$ & --      & $9.30031\mathrm{e}{+0}$ & --   & $6.61659\mathrm{e}{+1}$ & --      & $5.31510\mathrm{e}{-1}$ & --\\
0.17 & $3.69079\mathrm{e}{+3}$ & $-2.39$ & $6.09324\mathrm{e}{+0}$ & 0.60 & $2.06654\mathrm{e}{+2}$ & $-1.60$ & $1.64718\mathrm{e}{-1}$ & 1.65\\
0.08 & $1.21189\mathrm{e}{+3}$ & 1.54    & $3.23038\mathrm{e}{+0}$ & 0.88 & $3.87604\mathrm{e}{+1}$ & 2.31    & $4.18541\mathrm{e}{-2}$ & 1.89\\
0.04 & $3.78355\mathrm{e}{+2}$ & 1.66    & $1.61519\mathrm{e}{+0}$ & 0.99 & $5.94679\mathrm{e}{+0}$ & 2.67    & $1.03575\mathrm{e}{-2}$ & 1.99\\
0.02 & $1.54411\mathrm{e}{+2}$ & 1.28    & $8.18150\mathrm{e}{-1}$ & 0.97 & $1.22055\mathrm{e}{+0}$ & 2.27    & $2.88477\mathrm{e}{-3}$ & 1.83\\
\hline
\end{tabular}}
\end{table}

For $\kappa=1$, the two methods exhibit comparable numerical accuracy and essentially the same convergence behavior, see Table \ref{tab:kappa-1}. In particular, both methods achieve approximately first-order convergence in the $H^1$ seminorm and second-order convergence in the $L^2$ norm. This isotropic case provides a natural baseline for assessing the sensitivity of the two methods to increasing anisotropy.

To quantify the sensitivity of each method to anisotropy, we further introduce the anisotropy amplification factor
\begin{equation*}
	A_E^{\mathrm{M}}(\kappa)
	=
	\frac{E^{\mathrm{M}}(\kappa)}
	{E^{\mathrm{M}}(1)},
	\qquad
	\mathrm{M}\in\{\mathrm{VEM},\mathrm{SFVEM}\}.
\end{equation*}
On the finest mesh with $h=0.02$, the amplification factors corresponding to the two methods under the cases of $\kappa=10^2$ and $\kappa=10^4$ are
\begin{equation*}
\begin{aligned}
	A_{H^1}^{\mathrm{VEM}}(10^2)&=2.10,
	&\qquad
	A_{H^1}^{\mathrm{SFVEM}}(10^2)&=1.00,\\
	A_{H^1}^{\mathrm{VEM}}(10^4)&=189.74,
	&\qquad
	A_{H^1}^{\mathrm{SFVEM}}(10^4)&=1.01,
\end{aligned}
\end{equation*}
and
\begin{equation*}
\begin{aligned}
	A_{L^2}^{\mathrm{VEM}}(10^2)&=4.21,
	&\qquad
	A_{L^2}^{\mathrm{SFVEM}}(10^2)&=0.86,\\
	A_{L^2}^{\mathrm{VEM}}(10^4)&=435.89,
	&\qquad
	A_{L^2}^{\mathrm{SFVEM}}(10^4)&=0.90,
\end{aligned}
\end{equation*}
respectively.

The numerical results show a clear difference in the sensitivity of the two methods to anisotropy, see Tables \ref{tab:kappa-1}-\ref{tab:kappa-1e4}. As $\kappa$ increases, the errors of the standard VEM increase substantially. In particular, when $\kappa$ increases from $1$ to $10^4$, its $H^1$ and $L^2$ errors on the finest mesh are amplified by factors of approximately $190$ and $436$, respectively. By contrast, the corresponding amplification factors of stabilization-free VEM remain close to one, showing that its error are nearly unaffected by the increase of $\kappa$ in this experiment. 
Moreover, the stabilization-free VEM maintains approximately first-order convergence in the $H^1$ seminorm and second-order convergence in the $L^2$ norm for all tested values of $\kappa$. 
By contrast, for $\kappa=10^2$ and $10^4$, the standard VEM exhibits pronounced pre-asymptotic behavior on the coarse meshes: when the mesh is refined from $h=0.34$ to $0.17$, both errors increase, resulting in negative convergence rates. 
The errors of the standard VEM begin to decrease when $h\leq0.08$, but the observed rates remain affected by the anisotropy. 
These results demonstrate the substantially lower sensitivity of the stabilization-free VEM to strong anisotropy.

\subsection{Numerical comparison of stabilization-free VEM and standard FEM}

\begin{figure}[ht]
\centering
\includegraphics[width=0.5\textwidth]{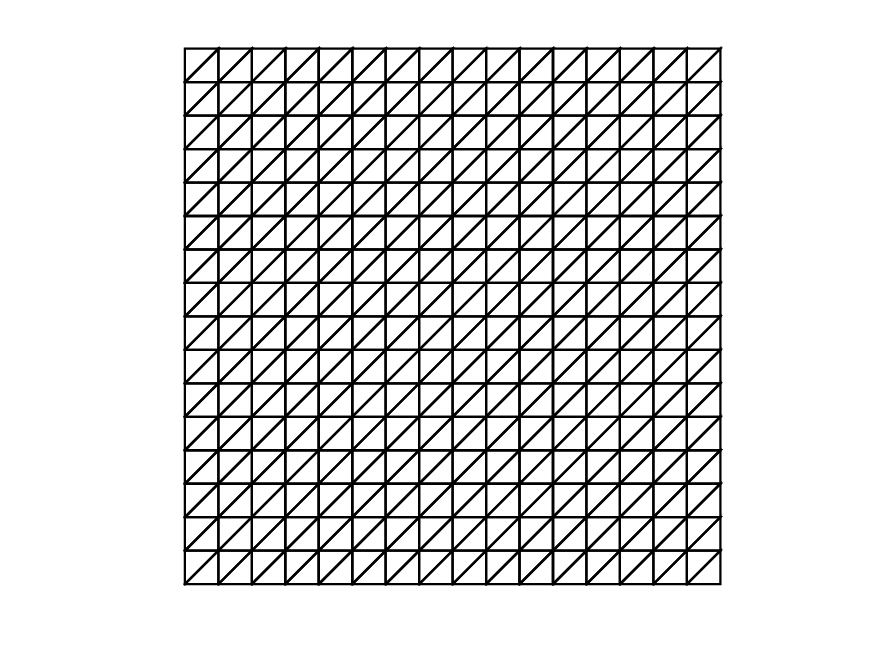}
\caption{A representative structured triangular mesh.}
\label{fig:structured-triangular-mesh}
\end{figure}

In this subsection, we compare the stabilization-free VEM with the standard FEM in terms of numerical accuracy and computational cost. Here we consider the second-order problem \eqref{eq:2ndEquation} with $\alpha=1$ and $\beta=0$. The right-hand side is chosen such that the exact solution is 
\[
u(x_1,x_2)=16x_1(1-x_1)x_2(1-x_2)\arctan\!\left(50\left(x_1-\tfrac12\right)\right),
\]
which exhibits a sharp vertical internal layer near $x_1=\tfrac12$. We solve this problem by using the stabilization-free conforming VEM and Lagrange FEM with the degree $k=2,3$ on five refined structured triangular meshes, with a representative mesh shown in Figure~\ref{fig:structured-triangular-mesh}. 

For a consistent comparison, the error of the stabilization-free VEM solution are taken to be $(u-\widetilde J_hu_h^{\mathrm{SFVEM}})$ with the reconstruction $\widetilde J_hu_h^{\mathrm{SFVEM}}$, since the VE solutions are not available in closed form. The reported computational time $T$ measures the total time required for the numerical discretization and the solution of the resulting linear system. The numerical results are reported in Tables~\ref{tab:poisson-critical-layer-beta0-k2} and~\ref{tab:poisson-critical-layer-beta0-k3}. 

\begin{table}[ht]
\centering
\caption{Numerical results for the internal-layer problem with $k=2$.}
\label{tab:poisson-critical-layer-beta0-k2}
\resizebox{\textwidth}{!}{%
\begin{tabular}{c|cc|cc|cc}
\hline
$h$ & $|u-u_h^{\mathrm{FEM}}|_1$ & $|u-\widetilde J_hu_h^{\mathrm{SFVEM}}|_1$ & $\lVert u-u_h^{\mathrm{FEM}}\rVert$ & $\lVert u-\widetilde J_hu_h^{\mathrm{SFVEM}}\rVert$ & $T^{\mathrm{FEM}}$ (s) & $T^{\mathrm{SFVEM}}$ (s) \\
\hline
0.3536  & $4.26794\mathrm{e}{+0}$ & $3.78732\mathrm{e}{+0}$ & $1.29433\mathrm{e}{-1}$ & $1.08756\mathrm{e}{-1}$ & 0.0026 & 0.0045 \\
0.1768  & $2.77175\mathrm{e}{+0}$ & $2.45984\mathrm{e}{+0}$ & $4.96299\mathrm{e}{-2}$ & $4.15452\mathrm{e}{-2}$ & 0.0035 & 0.0200 \\
0.08839 & $1.19099\mathrm{e}{+0}$ & $1.06142\mathrm{e}{+0}$ & $1.16321\mathrm{e}{-2}$ & $9.73856\mathrm{e}{-3}$ & 0.0115 & 0.0841 \\
0.04419 & $2.34696\mathrm{e}{-1}$ & $2.20545\mathrm{e}{-1}$ & $9.43543\mathrm{e}{-4}$ & $8.42077\mathrm{e}{-4}$ & 0.0496 & 0.3538 \\
0.02210 & $1.38169\mathrm{e}{-1}$ & $1.20208\mathrm{e}{-1}$ & $3.24281\mathrm{e}{-4}$ & $2.60696\mathrm{e}{-4}$ & 0.4969 & 2.1408 \\
\hline
\end{tabular}%
}
\end{table}

\begin{table}[ht]
\centering
\caption{Numerical results for the internal-layer problem with $k=3$.}
\label{tab:poisson-critical-layer-beta0-k3}
\resizebox{\textwidth}{!}{%
\begin{tabular}{c|cc|cc|cc}
\hline
$h$ & $|u-u_h^{\mathrm{FEM}}|_1$ & $|u-\widetilde J_hu_h^{\mathrm{SFVEM}}|_1$ & $\lVert u-u_h^{\mathrm{FEM}}\rVert$ & $\lVert u-\widetilde J_hu_h^{\mathrm{SFVEM}}\rVert$ & $T^{\mathrm{FEM}}$ (s) & $T^{\mathrm{SFVEM}}$ (s) \\
\hline
0.3536  & $2.60441\mathrm{e}{+0}$ & $2.30857\mathrm{e}{+0}$ & $5.88140\mathrm{e}{-2}$ & $4.87370\mathrm{e}{-2}$ & 0.0015 & 0.0119 \\
0.1768  & $1.10020\mathrm{e}{+0}$ & $9.59069\mathrm{e}{-1}$ & $1.42804\mathrm{e}{-2}$ & $1.17196\mathrm{e}{-2}$ & 0.0053 & 0.0340 \\
0.08839 & $2.52485\mathrm{e}{-1}$ & $2.48050\mathrm{e}{-1}$ & $1.19228\mathrm{e}{-3}$ & $1.14675\mathrm{e}{-3}$ & 0.0218 & 0.1415 \\
0.04419 & $1.76002\mathrm{e}{-1}$ & $1.53285\mathrm{e}{-1}$ & $5.53822\mathrm{e}{-4}$ & $4.58425\mathrm{e}{-4}$ & 0.0664 & 0.6243 \\
0.02210 & $3.82295\mathrm{e}{-2}$ & $3.25125\mathrm{e}{-2}$ & $6.27202\mathrm{e}{-5}$ & $5.07019\mathrm{e}{-5}$ & 0.8283 & 4.5931 \\
\hline
\end{tabular}%
}
\end{table}

The numerical results show that the FEM and stabilization-free VEM exhibit essentially the same convergence behavior for both $k=2$ and $k=3$. A more detailed observation reveals that the stabilization-free VEM yields smaller $H^1$- and $L^2$-errors than the standard FEM on each mesh with the same refinement level, which can be attributed to the larger number of DoFs in the stabilization-free VEM. Nevertheless, the stabilization-free VEM incurs higher computational cost than the FEM. The extra computational overhead of the VEM mainly stems from its local projection computations, finite element interpolation procedures, and its increased number of DoFs.

\section{Conclusion}

In this paper, we introduce a new framework for designing stabilization-free virtual element methods (VEMs) based on finite element interpolation with suitably chosen properties. Within this framework, we propose two types of stabilization-free schemes.
The first scheme employs a finite element interpolation that is computable, polynomial-preserving, norm-equivalent, and polynomial-consistent. As a result, it retains both consistency and stability in the same manner as standard VEMs.
The second scheme relies on an interpolation that is only required to be computable, polynomial-preserving, and norm-equivalent. This relaxation makes the scheme simpler to construct, requiring fewer degrees of freedom, and more computationally efficient in practice. Therefore, the second scheme is preferable in solving a broader range of problems, including those with nonlinearities or variable coefficients, where the consistency generally does not hold.
We derive abstract error estimates for both families of stabilization-free schemes. The theoretical framework is then used to develop conforming and nonconforming stabilization-free VEMs in both two and three dimensions. In contrast to existing gradient-projection-based approaches, the proposed FE-interpolation strategy simultaneously removes stabilization terms arising from both diffusion and reaction terms. Moreover, this new strategy is applicable to  other polytopal discretization settings, such as the hybrid high-order method and the weak Galerkin method.

The FE interpolations used in the framework are constructed for conforming and nonconforming VEMs in 2D and 3D. To realize the first type of stabilization-free VEMs (with consistency), we construct the polynomial-consistent FE interpolations in 2D and 3D. For the second type (without consistency), we build  the FE interpolations without polynomial consistency, again in both 2D and 3D. Finally, numerical experiments confirm that the proposed conforming and nonconforming stabilization-free VEMs achieve optimal convergence rates and exhibit significantly lower sensitivity to strong anisotropy compared to standard VEM.

\par\vspace{1.5em}
\noindent\textbf{Acknowledgements}\\[6pt]

 The work was supported by National Natural Science Foundation of China (Nos. 12371411, 12001170, 12271514, 12671486). The author Jikun Zhao was funded by Training Program for Young Backbone Teachers in Higher Education Institutions of Henan Province (No. 2025GGJS005). The author Bei Zhang was funded by Cultivation Program for Young Backbone Teachers in Henan University of Technology.

\end{document}